%% file: PreprintArxiv.tex
\documentclass[paper=a4,abstract=true,numbers=noenddot]{scrartcl}

\usepackage[english]{babel}  
\usepackage{hyperref}
\usepackage[intlimits]{amsmath}
\usepackage{amsthm}
\usepackage{amssymb} 
\usepackage{amsfonts,mathrsfs,bbm}
\allowdisplaybreaks[1]
\usepackage[numbers]{natbib}
\usepackage{graphicx}
\usepackage[space]{grffile} 
\usepackage{tikz}
\usepackage{pgfplots}
\usepackage[T1]{fontenc}
\usepackage{mathptmx}
\usepackage[scaled=.92]{helvet}
\usepackage{courier}
\usepackage[mathscr]{euscript}

\usepackage{authblk}
\usepackage{bm}
\usepackage{graphicx}
\usepackage{physics}
\usepackage{tensor}
\usepackage{xcolor}
\usepackage{enumerate}
\usepackage{algorithm}
\usepackage{algpseudocode}
\usepackage{subcaption}
\usepackage{siunitx}

\def\XXint#1#2#3{{\setbox0=\hbox{$#1{#2#3}{\int}$}
     \vcenter{\hbox{$#2#3$}}\kern-.5\wd0}}

\usetikzlibrary{patterns}
\pgfplotsset{compat=1.18}
\usepackage{array} 
\usepackage{booktabs} 
\usepackage{setspace}

\theoremstyle{remark}

\newcommand{\mc}[1]{\mathcal{#1}}
\newcommand{\mr}[1]{\mathrm{#1}}
\newcommand{\kl}[1]{\left(#1\right)}
\newcommand{\x}[0]{\mathbf{x}}
\newcommand{\vbf}[0]{\mathbf{v}}
\newcommand{\w}[0]{\mathbf{w}}
\newcommand{\ubf}[0]{\mathbf{u}}
\newcommand{\tbf}[0]{\mathbf{t}}
\newcommand{\g}[0]{\mathbf{g}}
\newcommand{\n}[0]{\mathbf{n}}
\newcommand{\y}[0]{\mathbf{y}}
\newcommand{\U}[0]{\mathbf{U}}
\newcommand{\Hbf}[0]{\mathbf{H}}
\newcommand{\f}[0]{\mathbf{f}}

\newcommand{\mat}[1]{\mathsf{#1}}   
\newcommand{\wmat}[1]{\widetilde{\mathsf{#1}}}
\newcommand{\hmat}[1]{\widehat{\mathsf{#1}}} 
\newcommand{\hwmat}[1]{\widehat{\widetilde{\mathsf{#1}}}}
\newcommand{\trce}[0]{\textup{Tr }} 
\newcommand{\dsy}[0]{\mathrm{d}s_\y}
\newcommand{\dsx}[0]{\mathrm{d}s_\x}  

\newcommand{\cvek}{\mat{c}}
\newcommand{\bvek}{\mat{b}}
\newcommand{\bone}{\mathbbm{1}} 
\newcommand{\A}{\mat{A}}

\newcommand{\bigO}{\mathcal{O}}

\newcommand{\etal}{\textup{et al. }}
\newcommand{\ie}{\textup{i.e., }}
\newcommand{\eg}{\textup{e.g., }}

\title{A 3D-ACA accelerated time domain boundary element method for elastodynamics using FMM and $\mathcal{H}$-matrix techniques}

\author[1]{V. Lakshmi Keshava}
\author[1]{M. Schanz}
\affil[1]{\mbox{Institute of Applied Mechanics, Graz University of Technology, Technikerstraße 4/II, 8010 Graz, Austria,} \protect\\
v.lakshmikeshava@tugraz.at, m.schanz@tugraz.at}
\date{}                     
\begin{document}
	
\maketitle
\section*{Abstract}
The time-domain Boundary Element Method (BEM) for linear elastodynamics with vanishing initial conditions is considered. Spatial discretization uses standard low-order boundary elements, while temporal discretization employs the generalized Convolution Quadrature (gCQ) method. The gCQ framework requires evaluating BEM matrices in the Laplace domain at several complex frequencies along a chosen contour, producing a three-dimensional tensor with one spatial matrix slice per frequency. To reduce storage and computational cost, a low rank approximation of the tensor is computed using 3D-Adaptive Cross Approximation (3D-ACA), extending the classical ACA to handle both the additional frequency dimension and the tensorial structure of elastodynamics. Within each frequency slice, the BEM matrices are further compressed using either the classical ACA algorithm using the $\mathcal{H}$-matrix approach or a Chebyshev interpolation based Fast Multipole Method (FMM).  A comparative study of all proposed methods is carried out using two academic examples, and the structural vibration of an induction machine is analyzed. 

\textbf{Keywords}: elastodynamics; boundary element method; generalised convolution quadrature; fast multipole method; multivariate adaptive cross approximation

\input{content}

\bibliography{library.bib}
\bibliographystyle{ieeetr}

\end{document}

%% file: content.tex
\section{Introduction}
\label{sec1}
The Boundary Element Method (BEM) is a widely used numerical tool for solving wave propagation problems in the time domain, particularly for semi-infinite domains and crack propagation problems. BEM is particularly advantageous for such geometries, as it reduces the dimensionality of the discretization and inherently satisfies the radiation condition in unbounded domains. The basis of BEM are the boundary integral equations with retarded potentials, derived through an analytical transformation of the underlying hyperbolic partial differential equation (PDE). This transformation relies on the so-called fundamental solution, which represents the exact solution of the PDE everywhere except at the origin, where it becomes singular. A comprehensive mathematical background on time-dependent boundary integral equations can be found in \cite{sayas16, costabel04}.

For elastodynamics, the first numerical realization using boundary integral formulation was published by Cruse and Rizzo \cite{cruse68}. They solved the transient problem in the Laplace transform space, utilizing the simpler fundamental solution of the corresponding elliptic PDE, and subsequently recovered the transient behaviour through numerical inversion. A similar Fourier domain transformation approach was presented by Dom\'{i}nguez \cite{dominguez78}, particularly for investigating the dynamic stiffness within an elastic half-space. These transformed domain formulations are memory efficient, since they rely only on the solution of the associated elliptic problems. Their main drawback, however, is that the complete time-dependent response can only be obtained by solving the transformed problem repeatedly over a large number of frequency points, which becomes computationally demanding, especially at high frequencies. The first boundary element formulation in time domain was developed by Mansur \cite{mansur83}, for zero initial conditions. This was later extended by Antes \cite{antes85} for non-zero initial conditions as well. In the early 1990s, Israil and Banerjee \cite{israil90} presented a more implementation-ready formulation, with higher order temporal interpolation and simplified kernels, with a focus on wave propagation in infinite and semi-infinite domains. In contrast to the former, these time-domain formulations require the explicit use of time-dependent fundamental solutions, leading to complicated kernel expressions and cumbersome numerical implementation. Furthermore, the memory requirement is also huge, since the spatial matrices have to be stored for each time-step. 

A good trade-off between these two approaches is provided by the BE formulations based on the convolution quadrature (CQ) method, introduced by Lubich \cite{lubich1, lubich2}. The essential idea is to discretize the temporal convolutions through the Laplace transform of the kernel combined with a linear multistep method, thereby requiring only the evaluation of the more simpler elliptic kernels. A good overview of the CQ based BEM for transient wave problems can be found in \cite{sayas16, schanz97, banjai12}. One major drawback with Lubich's CQ formulation is its restriction to uniform time stepping. This drawback was addressed by L\'{o}pez-Fern\'{a}ndez and Sauter \cite{lopez13, lopez15, lopez16}, who proposed a more generalized formulation called the generalized convolution quadrature method (gCQ), to handle variable time-step sizes. Applications of the gCQ-based time-domain BEM can be found in acoustics and thermoelasticity; see \eg \cite{sauter17, leitner21}. However, to the author's knowledge, there exist little to no such applications in elastodynamics. In this work, time-domain BE formulation based on gCQ method is used for elastodynamic applications. 

A major limitation of the BE formulations, regardless of the particular approach or problem class, is their quadratic complexity $\bigO \kl{M^2}$, in both storage and computation time for $M$ unknowns. In the case of gCQ-method based time-domain BEM, the formulation essentially gives rise to a three-dimensional data array that needs to be computed and stored. This structure arises because the two-dimensional array associated with the spatial discretization must be evaluated at each selected complex frequency. Consequently, the additional frequency dimension in the CQ-based formulation leads to an overall complexity of order $\bigO \kl{M^2N}$, where $N$ denotes the number of complex frequencies. For the spatial discretization, \ie the elliptic problem, several classes of fast methods have been developed over the last three decades. Among these, the more widely known are the kernel-dependent methods, such as the fast multipole method (FMM) \cite{greengard87, liu09, of05} and the pre-corrected fast Fourier transform \cite{phillips02, masters04} as well as the kernel-independent methods such as the SVD-based \cite{kapur97} and QR-based \cite{gope05} compression methods, $\mc{H}$-matrix based adaptive cross approximation (ACA) methods \cite{bebendorf03, bebendorf08, haider19}, to name only a few. 

On the contrary, comparatively fewer fast methods have been proposed or tested for the frequency dimension of the three-dimensional data array arising in the CQ-based formulations. There are, however, examples in which fast methods such as ACA \cite{messner10} or FMM \cite{maruyama16} have been employed to accelerate and compress the spatial dimension of CQ-based BEM. It is also worth mentioning the work of Takahashi et al. \cite{takahashi04}, who developed a fast method based on plane wave time domain (PWTD) algorithm for 3D-elastodynamics in time domain, which was an extension of the original idea proposed by Michielssen's group \cite{ergin98} for the scalar wave equation.

In the present work, a combined approach is adopted to achieve true speed up and compression in all dimensions. The frequency dimension is compressed using the so called multivariate adaptive cross approximation (3D-ACA), also known as generalized adaptive cross approximation. This method is a generalized extension of the original ACA, developed by Bebendorf et al. \cite{bebendorf11, bebendorf13}. For the spatial dimension, either the $\mc{H}$-matrix based ACA or the Chebyshev interpolation based black-box FMM \cite{fong09} is employed, and a comparative study of the two approaches is carried out. A similar study was performed by Schanz et al. \cite{schanz26} for the scalar wave equation. The purpose of this paper is to extend this framework to the matrix-valued elastodynamics in 3D. As will be seen, this extension is not straightforward, since the 3D-ACA in the frequency dimension has to be adapted to the matrix-valued elastodynamics, based on the work of Rjasanow and Weggler \cite{rjasanow17}. The use of this matrix-valued ACA can already be seen in the works of Chaillat \etal \cite{chaillat17} and Haider \etal \cite{haider19}. The work of Seibel \cite{seibel22} should also be mentioned here, where a similar approach was used for the CQ-based BEM applied to the scalar wave equation. In contrast to the present work, however, the spatial dimension there is approximated using the so-called $\mc{H}^2$ technique.  

In this paper, both collocation and Galerkin formulations are employed to demonstrate the advantages of the 3D-ACA-accelerated, gCQ-based time-domain BEM, first on some simple academic examples. A thorough comparative study is then carried out between the $\mc{H}$-matrix approach and the FMM approaches, and the investigation is further extended to a real-world application, \ie to compute the vibrations in an electric machine.


\section{Boundary element formulation}
\label{sec2}

\subsection{Problem statement}
\label{sec2.1}

Consider a bounded Lipschitz domain $\Omega \subset \mathbb{R}^3$ with its boundary denoted by $\Gamma=\partial\Omega$. A linear isotropic elastodynamics problem with no body forces is considered, with the displacement field denoted by $\ubf\kl{\x}$ and the traction vector denoted by $\tbf\kl{\x}$. The governing equations for the considered mixed initial boundary value problem can be written as
\begin{subequations}\label{eq:gov}
\begin{align}
-\kl{c_1^2-c_2^2} \nabla \kl{\nabla\cdot \ubf\kl{\x,t}} - c_2^2 \Delta \ubf\kl{\x,t} + \frac{\partial^2\ubf}{\partial t^2} \kl{\x,t} &= \mathbf{0} \qquad\qquad (\x,t) \in \Omega \times (0,T)\\
\trce\ubf\kl{\x,t} &= \g_D\kl{\x,t} \;\quad \kl{\x,t} \in \Gamma_D \times (0,T)\\
\kl{\mc{T}_\x\ubf}\kl{\x,t} &= \g_N\kl{\x,t} \;\quad \kl{\x,t} \in \Gamma_N \times (0,T)\\
\ubf\kl{\x,0} = \frac{\partial\ubf}{\partial t}\kl{\x,0} &= \mathbf{0}  \qquad\qquad \: \x\in \Omega.
\end{align}
\end{subequations}
The quantities $c_1$ and $c_2$ denote the longitudinal and shear wave speeds, respectively, and are given by 
\begin{equation}\label{eq:lame}
c_1 = \sqrt{\frac{\lambda+2\mu}{\rho}}, \quad c_2=\sqrt{\frac{\mu}{\rho}},
\end{equation}
where $\rho$ is the material density, and $\lambda$ and $\mu$ are the Lam\'{e} constants. The respective traces are defined by
\begin{equation}\label{eq:trace1}
\trce_\x\ubf\kl{\x,t} := \lim_{\Omega \ni  \x \to \x \in \Gamma}\ubf\kl{\x,t}
\end{equation}
for the Dirichlet trace and by
\begin{equation}\label{eq:traction}
\kl{\mc{T}_\x\ubf}\kl{\x,t}:= \lim_{\Omega \ni  \x \to \x \in \Gamma}[\bm{\sigma}\kl{\x,t}\cdot \n\kl{\x}] = \tbf\kl{\x,t}
\end{equation} 
for the Neumann trace, also known as the elastic traction operator, with the stress tensor $\bm{\sigma}\kl{\x,t}$ and the outer normal vector $\n\kl{\x}$ on the boundary $\Gamma$. The Cauchy stress tensor $\bm{\sigma}\kl{\x,t}$ is defined as
\begin{equation}
    \bm{\sigma}\kl{\x,t} = \rho\kl{c_1^2-2c_2^2}\kl{\nabla\cdot\ubf\kl{\x,t}}\mathbf{I} + \rho c_2^2 \kl{ \nabla\ubf\kl{\x,t}+\kl{\nabla\ubf\kl{\x,t}}^\top },
\end{equation}
$\mathbf{I}$ being the identity matrix. The Lipschitz boundary $\Gamma$ is split into two mutually disjoint sets $\Gamma_D$ and $\Gamma_N$, such that $\Gamma=\Gamma_D \cup \Gamma_N$. The surface displacements and tractions are prescribed by the given data $\g_D\kl{\x,t}$ on $\Gamma_D$ and $\g_N\kl{\x,t}$ on $\Gamma_N$ respectively, over the time interval $\kl{0,T}$, with $T>0$. For vector-valued elastodynamic problems, the boundary conditions may, in general, be prescribed component-wise, so that different boundary condition types can be assigned in different spatial directions at the same boundary point, as in the case of a roller support. For simplicity, such component-wise mixed boundary conditions are not considered in the present work.

\subsection{Boundary integral equation}
\label{sec2.2}

The problem given in \eqref{eq:gov} can be solved by means of boundary integral equations \cite{cruse68, dominguez93, bonnet99}, using either the indirect or the direct formulations, but only the latter is considered here. Based on the \textit{dynamic reciprocity theorem} \cite{graffi54, wheeler68}, the dynamic representation formula for the problem \eqref{eq:gov} can be written as
\begin{equation}\label{eq:rep}
\begin{aligned}
\ubf\kl{\x,t} &= \int_0^t\int_\Gamma \U^*\kl{\y-\x,t-\tau} \tbf\kl{\y,\tau} \dsy \mathrm{d}\tau \\
&- \int_0^t\int_\Gamma [\kl{\mc{T}_\y\U^*}\kl{\y-\x,t-\tau}]^\top \ubf\kl{\y,\tau}\dsy \mathrm{d}\tau \quad (\x,t) \in \Omega \times (0,T),
\end{aligned}
\end{equation}
where $\U^*$ is the time domain displacement fundamental solution. Applying the Dirichlet trace $\trce_\x$ from \eqref{eq:trace1} to the representation formula yields the standard boundary integral equation (SBIE),
\begin{equation}\label{eq:SBIE}
\mc{C}\kl{\x}\ubf\kl{\x,t} = \kl{\mc{V}\ast\tbf}\kl{\x,t} - \kl{\mc{K}\ast\ubf}\kl{\x,t} \quad (\x,t) \in \Gamma \times (0,T)
\end{equation}
where the limiting process gives rise to $\mc{C}\kl{\x}$, the so-called integral-free term which depends on the local geometry (the solid angle) at $\x$ in collocation BEM \cite{mantic93}. In the Galerkin formulation, however, this term simplifies to $\mc{C}\kl{\x}=\frac{1}{2}\mathbf{I}$. A symmetric Galerkin formulation additionally requires the hypersingular boundary integral equation (HBIE), which is obtained by applying the traction operator \eqref{eq:traction} to the representation formula \eqref{eq:rep}. The HBIE is given by 
\begin{equation}\label{eq:HBIE}
[\mathbf{I}-\mc{C}\kl{\x}]\tbf\kl{\x,t} = \kl{\mc{K}'\ast\tbf}\kl{\x,t} + \kl{\mc{D}\ast\ubf}\kl{\x,t} \quad (\x,t) \in \Gamma \times (0,T).
\end{equation}  
Corresponding to these integral equations, the associated boundary integral operators are introduced. The single-layer potential is defined by
\begin{equation}\label{eq:SLP}
\kl{\mc{V}\ast\tbf}\kl{\x,t} := \int_0^t\int_\Gamma \U^*\kl{\y-\x,t-\tau} \tbf\kl{\y,\tau} \dsy \mathrm{d}\tau.
\end{equation}
The double-layer and the adjoint double-layer potentials are given by
\begin{equation}\label{eq:DLP}
\kl{\mc{K}\ast\ubf}\kl{\x,t} := \int_0^t\int_\Gamma [\kl{\mc{T}_\y\U^*}\kl{\y-\x,t-\tau}]^\top \ubf\kl{\y,\tau}\dsy \mathrm{d}\tau,
\end{equation}
\begin{equation}\label{eq:ADLP}
\kl{\mc{K'}\ast\tbf}\kl{\x,t} := \int_0^t \mc{T}_\x \int_\Gamma \U^*\kl{\y-\x,t-\tau} \tbf\kl{\y,\tau}\dsy \mathrm{d}\tau.
\end{equation}
Finally, the hypersingular operator is defined as
\begin{equation}\label{eq:HSO}
\kl{\mc{D}\ast\ubf}\kl{\x,t} := - \int_0^t \mc{T}_\x \int_\Gamma [\kl{\mc{T}_\y\U^*}\kl{\y-\x,t-\tau}]^\top \ubf\kl{\y,\tau}\dsy \mathrm{d}\tau.
\end{equation}
The single-layer potential is weakly singular, whereas the double-layer and the adjoint double-layer potentials are strongly singular due to the application of the traction operator on $\U^*$, and must be understood in the sense of a Cauchy principal value. Lastly, as the name suggests, the hypersingular operator is hypersingular, and therefore be interpreted as a finite part integral, in the sense of Hadamard \cite{hadamard23}. Both the strongly singular and hypersingular integrals are regularized to a weakly singular form by means of partial integration \cite{kielhorn08}. 

\subsection{Spatial discretization}
\label{sec2.3}  
The boundary $\Gamma$ is discretized using flat triangular elements, resulting in an approximate boundary $\Gamma^h$, given by
\begin{equation}\label{eq:tri}
\Gamma \approx \Gamma^h = \bigcup_{k=1}^{K}\gamma_k.
\end{equation}
Thus, the approximate boundary $\Gamma^h$ is a union of $K$ triangular elements $\gamma_k$, on which the unknown boundary data are associated with either $\Gamma_D$ or $\Gamma_N$. The mixed boundary conditions in \eqref{eq:gov} imply that the boundary fields $\ubf\kl{\x,t}$ and $\tbf\kl{\x,t}$ are partially unknown, depending on the boundary segment. Accordingly, they are decomposed as 
\begin{equation}\label{eq:decomp}
\ubf\kl{\x,t} = \ubf^h\kl{\x,t} + \g^h_D\kl{\x,t} \quad \textup{and} \quad \tbf\kl{\x,t} = \tbf^h\kl{\x,t} + \g^h_N\kl{\x,t} \quad \textup{for } \x \in \Gamma,
\end{equation}
where $\g^h_D$ and $\g^h_N$ are arbitrary but fixed extensions of the prescribed boundary data, chosen to vanish on the complementary boundary segments. The corresponding spaces are defined by
\begin{equation}\label{eq:subspace}
\begin{aligned}
S_h\kl{\Gamma_D}&=\textup{span}\{\varphi_1, \varphi_2,...,\varphi_{M_1}\}\\
S_h\kl{\Gamma_N}&=\textup{span}\{\psi_1, \psi_2,...,\psi_{M_2}\},
\end{aligned}
\end{equation} 
where $M_1$ continuous, piecewise linear shape functions $\varphi$ are used for the approximation of the Dirichlet unknowns, and $M_2$ discontinuous, piecewise constant shape functions $\psi$ are used for the approximation of Neumann unknowns. The corresponding approximations are then given by 
\begin{equation}\label{eq:basis}
\ubf^h\kl{\x,t} = \sum_{l=1}^{M_1} \ubf_l\kl{t} \varphi_l\kl{\x} \quad \textup{and} \quad \tbf^h\kl{\x,t} = \sum_{j=1}^{M_2} \tbf_j\kl{t} \psi_l\kl{\x}.
\end{equation}
The coefficients $\ubf_l\kl{t}$ and $\tbf_j\kl{t}$ are still continuous in time and, therefore, still require temporal discretization. Inserting the above shape functions into the SBIE \eqref{eq:SBIE}, the semi-discretized equation can be written in the matrix form as
\begin{equation}\label{eq:discreteCollo}
\begin{bmatrix}
\mat{V}_{D}\kl{t} & -\mat{K}_{D}\kl{t} \\
\mat{V}_{N}\kl{t} & -\kl{\mat{C}_{N}+\mat{K}_{N}\kl{t}}
\end{bmatrix} \ast 
\begin{bmatrix}
\mat{t}^h_D\kl{t} \\
\mat{u}^h_N\kl{t}
\end{bmatrix} = 
\begin{bmatrix}
\mat{C}_{D}+\mat{K}_{D}\kl{t} & -\mat{V}_{D}\kl{t} \\
\mat{K}_{N}\kl{t} & -\mat{V}_{N}\kl{t}
\end{bmatrix} \ast 
\begin{bmatrix}
\mat{g}^h_D\kl{t} \\
\mat{g}^h_N\kl{t}
\end{bmatrix},
\end{equation}
where $\mat{V}$ and $\mat{K}$ denote the semi-discrete single-layer and double-layer collocation matrices, respectively, given by 
\begin{equation}\label{discrVK}
\begin{aligned}
\mat{V}\kl{t}[i,j] &= \int_{supp\kl{\psi_j}}\U^*\kl{\y_i -\x, t}\psi_j\kl{\y}\dsy \\
\mat{K}\kl{t}[i,j] &= \int_{supp\kl{\varphi_j}}[\kl{\mc{T}_\y\U^*}\kl{\y_i -\x, t}]^\top \varphi_j\kl{\y}\dsy.
\end{aligned}
\end{equation}
Here, the subscripts \(D\) and \(N\) indicate the corresponding collocation nodes at $\Gamma_D$ and $\Gamma_N$, respectively. The collocation points $\x_i$ are chosen as the nodes of the linear triangular elements for $\x_i \in \Gamma_N$, and as the midpoints of the constant triangular elements for $\x_i \in \Gamma_D$. 

For the symmetric Galerkin formulation, both the SBIE \eqref{eq:SBIE} and the HBIE \eqref{eq:HBIE} are tested with appropriate test functions $\vbf\kl{\x}$ and $\w\kl{\x}$, which leads to the variational formulation 
\begin{equation}\label{eq:var}
\begin{aligned}
\langle \mc{V} \ast \mat{t}^h, \w \rangle_{\Gamma_D} - \langle \mc{K} \ast \mat{u}^h, \w \rangle_{\Gamma_D} &= \langle (\frac{1}{2}\mc{I}+\mc{K}) \ast \mat{g}^h_D - \mc{V} \ast \mat{g}^h_N, \w \rangle_{\Gamma_D} \\ 
\langle \mc{K}' \ast \mat{t}^h, \vbf \rangle_{\Gamma_N} + \langle \mc{D} \ast \mat{u}^h, \vbf \rangle_{\Gamma_N} &= \langle (\frac{1}{2}\mc{I}-\mc{K}') \ast \mat{g}^h_N - \mc{D} \ast \mat{g}^h_D, \vbf \rangle_{\Gamma_N}.
\end{aligned}
\end{equation}
Here, $\mc{I}$ denotes the identity operator 
\begin{equation}\label{eq:identity}
\kl{\mc{I}\ast \ubf}_{\Gamma}\kl{\x,t} = \int_0^t \int_\Gamma \delta\kl{\y-\x,t-\tau}\mathbf{I}\cdot \ubf\kl{\y,\tau} \dsy \mathrm{d}\tau
\end{equation}
and the inner product $\langle f,g \rangle_\Gamma = \int_\Gamma f(x) g(x) \dsx$ implies that the Galerkin formulation involves double spatial integrations. Note that the test functions $\vbf\kl{\x}$ and $\w\kl{\x}$ do not depend on time and are approximated by the same polynomial spaces as those used for the Dirichlet and Neumann data in \eqref{eq:basis}.

Inserting the spatial approximations \eqref{eq:basis} into the \eqref{eq:var} results in the semi-discrete equation system 
\begin{equation}\label{eq:discreteGal}
\begin{bmatrix}
\wmat{V}_{D}\kl{t} & -\wmat{K}_{D}\kl{t} \\
\wmat{K}_{N}'\kl{t} & \wmat{D}_{N}\kl{t}
\end{bmatrix} \ast 
\begin{bmatrix}
\mat{t}^h_D\kl{t} \\
\mat{u}^h_N\kl{t}
\end{bmatrix} = 
\begin{bmatrix}
(\frac{1}{2}\widetilde{\mathscr{I}}_{D}\kl{t}+\wmat{K}_{D}\kl{t}) & -\wmat{V}_{D}\kl{t} \\
(\frac{1}{2}\widetilde{\mathscr{I}}_{N}\kl{t}-\wmat{K}_{N}'\kl{t}) & -\wmat{D}_{N}\kl{t}
\end{bmatrix} \ast 
\begin{bmatrix}
\mat{g}^h_D\kl{t} \\
\mat{g}^h_N\kl{t}
\end{bmatrix}.
\end{equation}
Here, $\wmat{V}$ denotes the semi-discrete single-layer Galerkin matrix, given by
\begin{equation}\label{eq:GalV}
\wmat{V}\kl{t}[i,j] = \int_{supp\kl{\psi_i}}\psi_i\kl{\x}\int_{supp\kl{\psi_j}}\U^*\kl{\y -\x, t}\psi_j\kl{\y}\dsy\dsx.
\end{equation}
Similarly, $\wmat{K}$ and $\wmat{K}'$ denote the semi-discrete double-layer and adjoint double-layer Galerkin matrices
\begin{subequations}\label{eq:GalKK}
\begin{align}
\wmat{K}\kl{t}[i,j] &= \int_{supp\kl{\psi_i}}\psi_i\kl{\x}\int_{supp\kl{\varphi_j}}[\kl{\mc{T}_\y\U^*}\kl{\y -\x, t}]^\top \varphi_j\kl{\y}\dsy\dsx \\
\wmat{K}'\kl{t}[i,j]&= \int_{supp\kl{\varphi_i}} \varphi_i\kl{\x} \mc{T}_\x \int_{supp\kl{\psi_j}}\U^*\kl{\y -\x, t}\psi_j\kl{\y}\dsy\dsx,
\end{align}
\end{subequations}
and $\wmat{D}$ denotes the semi-discrete hypersingular Galerkin matrix, 
\begin{equation}\label{eq:GalD}
\wmat{D}\kl{t}[i,j] = - \int_{supp\kl{\varphi_i}} \varphi_i\kl{\x} \mc{T}_\x \int_{supp\kl{\varphi_j}}[\kl{\mc{T}_\y\U^*}\kl{\y -\x, t}]^\top\varphi_j\kl{\y}\dsy\dsx.
\end{equation}
Lastly, the entries of the identity operator are computed as 
\begin{equation}\label{eq:Identity}
\widetilde{\mathscr{I}}[i,j] = \int_{supp\kl{\varphi_i}} \varphi_i\kl{\x}\psi_j\kl{\x}\dsx.
\end{equation}
The strongly singular and hypersingular integrals are treated analogously to the collocation method. Weakly singular integrals in the Galerkin formulation are evaluated numerically using the approach by Erichsen and Sauter \cite{erichsen98}, while in the collocation method they are treated by means of the Duffy transformation \cite{duffy82}. 

\subsection{Temporal discretization}
\label{sec2.4} 

The semi-discrete integral equations are discretized in time by the Runge--Kutta based generalized convolution quadrature (gCQ) method \cite{lopez13, lopez16}. Considering a single entry of the matrix--vector convolution appearing in the semi-discrete system \eqref{eq:discreteCollo}, the convolution integral can be written as
\begin{equation}\label{eq:matvec}
\kl{\mat{V}_{D}[i,j] \ast \mat{t}_D^h[j]}\kl{t} = \int_0^t \mat{V}_D[i,j]\kl{t-\tau} \mat{t}_D^h[j]\kl{\tau}\mr{d}\tau = \mat{f}\kl{t}.
\end{equation}
Using the inverse Laplace transform, this convolution integral can be reformulated as
\begin{equation}\label{eq:laptrans}
\begin{aligned}
\mat{f}\kl{t} & = \int_0^t \mat{V}_D[i,j]\kl{t-\tau} \mat{t}_D^h[j]\kl{\tau}\mr{d}\tau \\
&= \frac{1}{2\pi\mr{i}}\int_C \hmat{V}_{D}[i,j]\kl{s} \int_0^t e^{s\kl{t-\tau}}\mat{t}_D^h[j]\kl{\tau}  \mr{d}\tau \mr{ds}, 
\end{aligned}
\end{equation} 
where $s$ denotes the Laplace variable with $s \in \mathbb{C}, s.t. \mathfrak{R}s>0$. The notation $\hmat{V}_{D}[i,j]$ indicates the Laplace transform of $\mat{V}_{D}[i,j]$, which means that the temporal fundamental solution $\U^*\kl{\y-\x, t-\tau}$ is replaced by its Laplace-domain counterpart $\U\kl{\y-\x, s}$. The outer integral is a contour integral taken along the Bromwich contour $C$, oriented from $a-\mr{i}\infty \to a+\mr{i}\infty$, where $a$ is chosen to lie to the right of all singularities of $\hmat{V}_{D}[i,j]$. The inner integral can then be defined as 
\begin{equation}\label{eq:ODE1}
x\kl{s,t} := \int_0^t e^{s\kl{t-\tau}}\mat{t}_D^h[j]\kl{\tau}  \mr{d}\tau,
\end{equation}
where $x$ satisfies the ordinary differential equation
\begin{equation}\label{eq:ODE2}
\frac{\partial x\kl{s,t}}{\partial t} = sx\kl{s,t} + \mat{t}_D^h[j]\kl{t}, \quad x\kl{s,0}=0.
\end{equation}
This ordinary differential equation is then solved numerically by a time-stepping scheme with not necessarily constant time steps. To this end, an A- and L-stable Runge--Kutta method of $m$ stages, characterized by the Butcher table $\A= \kl{a_{i,j}}^{m}_{i,j=1}$, $\bvek=\kl{b_{i}}^{m}_{i=1}$, $\cvek=\kl{c_{i}}^{m}_{i=1}$ is considered. For the temporal discretization, a sequence of time points 
\begin{equation*}
[0,T] =[0,t_1,t_2,...,t_N], \quad \Delta t_j = t_j - t_{j-1}, \quad j=1,2,...,N.,
\end{equation*}  
is employed. The stability assumptions require that the stability function, defined by 
\begin{equation}\label{eq:stability}
R\kl{z} := 1 + z\bvek^\top \kl{\mathbf{I}-z\A}^{-1}\bone, \quad \bone:= \kl{1,1...,1}^\top
\end{equation}
is bounded as 
\begin{equation}\label{eq:bounded}
|R\kl{z}| \leq 1, \quad \textup{and} \quad \kl{\mathbf{I}-z\A} \textup{is nonsingular for all } \mathfrak{R}z\leq 0,
\end{equation}
and that $\bvek^\top\A^{-1}=\kl{0,0,...1}$ holds. For further details on these restrictions and assumptions, the reader is referred to \cite{lopez16, schanz01}. With these definitions, the gCQ algorithm for the matrix-vector product in \eqref{eq:matvec} can be written as in \cite{lopez15}:
\begin{itemize}
\item Step $n=1$
\begin{equation*}
\mat{f}\kl{t}= \hmat{V}_{D}\kl{\kl{\Delta t_{n=1} \mat{A}}^{-1}}(\mat{t}_D^h)_{n=1},
\end{equation*}
where zero initial conditions are assumed implicitly.
\item Steps $n=2,3,...,N$ 

\begin{enumerate}
\item Update the solution vector $\mat{x}_{n-1}$ at all integration points $s_{\ell}$:
\begin{equation}\label{eq:ODEsol}
\mat{x}_{n-1}\kl{s_\ell} = \kl{\mathbf{I}-\Delta t_{n-1}s_\ell\A}^{-1}  \kl{\kl{\bvek^\top\A^{-1}\cdot \mat{x}_{n-2}\kl{s_\ell}}\bone + \Delta t_{n-1}\A (\mat{t}_D^h)_{n-1}}   
\end{equation}
for $\ell=1,...,N_Q$ with $N_Q$ being the number of integration points. 

\item Solve the following system at the current time step $t_n$
\begin{equation}\label{eq:gCQ}
\kl{\mat{f}}_n = \hmat{V}_{D}\kl{\kl{\Delta t_{n} \mat{A}}^{-1}}(\mat{t}_D^h)_{n} + \sum_{\ell=1}^{N_Q}\hmat{V}_{D}\kl{s_\ell}\mat{W}_\ell^{\Delta t_n}\kl{(\mat{t}_D^h)_{n-1}}   
\end{equation}
\end{enumerate}
\end{itemize}

Equation \eqref{eq:gCQ} gives the computation of one entry in \eqref{eq:matvec} at all $m$ stages, where $(\mat{t}_D^h)_{n}$ denotes the vector of all nodal values of $\mat{t}_D^h$ at time $t_n$, collected over all $m$ stages. Furthermore, the notation $\hmat{V}_{D}\kl{\kl{\Delta t_{n} \mat{A}}^{-1}}$ indicates that the discrete operator $\hmat{V}_{D}$ is evaluated at all values of ${\kl{\Delta t_{n} \mat{A}}^{-1}}$, so that each entry becomes of size $m\times m$. Finally, the term $$\mat{W}_\ell^{\Delta t_n}\kl{(\mat{t}_D^h)_{n-1}} = \omega_\ell \kl{\bvek^\top\A^{-1}\cdot \mat{x}_{n-1}\kl{s_\ell}}  \kl{\mathbf{I}-\Delta t_{n}s_\ell\A}^{-1}  \bone$$ is a vector of size $M_2m$, where $\omega_\ell$ are the integration weights and $\mat{x}_{n-1}\kl{s_\ell}$ from \eqref{eq:ODEsol} is again a vector collecting the solutions at all $m$ stages for each node in $(\mat{t}_D^h)$.

With the above temporal discretization, the discrete sets of integral equations can be written explicitly for the collocation approach as 
\begin{equation}
\begin{aligned}\label{eq:discrCollo}
\begin{bmatrix}
\hmat{V}_{D} & -\hmat{K}_{D} \\
\hmat{V}_{N} & -\kl{\mat{C}_{N}+\hmat{K}_{N}}
\end{bmatrix} &\kl{\kl{\Delta t_{n} \mat{A}}^{-1}}
\begin{bmatrix}
\mat{t}^h_D \\
\mat{u}^h_N
\end{bmatrix}_n = 
\begin{bmatrix}
\mat{C}_{D}+\hmat{K}_{D} & -\hmat{V}_{D} \\
\hmat{K}_{N} & -\hmat{V}_{N}
\end{bmatrix} \kl{\kl{\Delta t_{n} \mat{A}}^{-1}}
\begin{bmatrix}
\mat{g}^h_D \\
\mat{g}^h_N
\end{bmatrix}_n     \\
 & \: \hspace{-30mm}
 + \sum_{\ell=1}^{N_Q}\kl{\begin{bmatrix}
\hmat{K}_{D} & -\hmat{V}_{D} \\
\hmat{K}_{N} & -\hmat{V}_{N}
\end{bmatrix}  \kl{s_\ell} 
\begin{bmatrix}
\mat{W}_\ell^{\Delta t_n}\kl{(\mat{g}_D^h)_{n-1}}  \\
\mat{W}_\ell^{\Delta t_n}\kl{(\mat{g}_N^h)_{n-1}} 
\end{bmatrix} - 
\begin{bmatrix}
\hmat{V}_{D} & -\hmat{K}_{D} \\
\hmat{V}_{N} & -\hmat{K}_{N}
\end{bmatrix}  \kl{s_\ell} 
\begin{bmatrix}
\mat{W}_\ell^{\Delta t_n}\kl{(\mat{t}_D^h)_{n-1}}  \\
\mat{W}_\ell^{\Delta t_n}\kl{(\mat{u}_N^h)_{n-1}} 
\end{bmatrix}}.
\end{aligned}    
\end{equation}
Similarly, the discrete Galerkin formulation can be written as 
\begin{equation}\label{eq:discrGal}
\begin{aligned}
\begin{bmatrix}
\hwmat{V}_{D} & -\hwmat{K}_{D} \\
\hwmat{K}_{N}' & \hwmat{D}_{N}
\end{bmatrix}  &\kl{\kl{\Delta t_{n} \mat{A}}^{-1}}
\begin{bmatrix}
\mat{t}^h_D \\
\mat{u}^h_N
\end{bmatrix}_n = 
\begin{bmatrix}
(\frac{1}{2}\widetilde{\mathscr{I}}_{D}+\hwmat{K}_{D}) & -\hwmat{V}_{D} \\
(\frac{1}{2}\widetilde{\mathscr{I}}_{N}-\hwmat{K}_{N}') & -\hwmat{D}_{N}
\end{bmatrix} \kl{\kl{\Delta t_{n} \mat{A}}^{-1}}
\begin{bmatrix}
\mat{g}^h_D \\
\mat{g}^h_N
\end{bmatrix}_n     \\
 & \: \hspace{-30mm}
+ \sum_{\ell=1}^{N_Q}\kl{\begin{bmatrix}
\hwmat{K}_{D} & -\hwmat{V}_{D} \\
-\hwmat{K}_{N}' & -\hwmat{D}_{N}
\end{bmatrix}  \kl{s_\ell} 
\begin{bmatrix}
\mat{W}_\ell^{\Delta t_n}\kl{(\mat{g}_D^h)_{n-1}}  \\
\mat{W}_\ell^{\Delta t_n}\kl{(\mat{g}_N^h)_{n-1}} 
\end{bmatrix} - 
\begin{bmatrix}
\hwmat{V}_{D} & -\hwmat{K}_{D} \\
\hwmat{K}_{N}' & \hwmat{D}_{N}
\end{bmatrix}  \kl{s_\ell} 
\begin{bmatrix}
\mat{W}_\ell^{\Delta t_n}\kl{(\mat{t}_D^h)_{n-1}}  \\
\mat{W}_\ell^{\Delta t_n}\kl{(\mat{u}_N^h)_{n-1}} 
\end{bmatrix}}.
\end{aligned}    
\end{equation}
The resulting system of equations requires the evaluation and storage of the Laplace-transformed integral kernels at $N_Q$ number of complex frequencies $s_\ell$, which gives rise to an array of system matrices of size $M\times M \times N_Q$, where $M$ denotes either $M_1$ or $M_2$ from \eqref{eq:subspace}, depending on the corresponding integral kernel. Each such matrix array in \eqref{eq:discrCollo} and \eqref{eq:discrGal} has a computational complexity of $\bigO\kl{\kl{N_Q+1}M^2}$, while the evaluation of the time-stepping scheme requires $\bigO\kl{N_QN}$ matrix-vector multiplications. These matrix arrays are therefore approximated by a data-sparse representation based on multivariate ACA (3D-ACA) \cite{bebendorf13}. In addition, each spatial matrix block in the resulting three-dimensional array is approximated either by a $\mc{H}$-matrix based classical ACA or by a Chebyshev-interpolation-based FMM.


\section{Multivariate adaptive cross approximation}
\label{sec3}

The 3D-ACA algorithm used in the present work was originally proposed in \cite{bebendorf13} for the approximation of a three-dimensional data array $C \in \mathbb{C}^{M\times M \times N_Q}$, and can be understood as an extension of the classical adaptive cross approximation. In order to account for the tensor-valued nature of the elastodynamic problem considered here, the matrix-valued ACA proposed in \cite{rjasanow17} is incorporated into this framework. 

The data-sparse approximation of the three-dimensional array $C$ is sought in the form of an outer product sum
\begin{equation}
C \approx C^r = \sum_{d=1}^{r} \Hbf_d \otimes \f_d,   
\end{equation}
where each $\Hbf_d \in \mathbb{C}^{M\times M}$ is a matrix and each $\f_d \in \mathbb{C}^{N_Q}$ is a vector. In the present setting, the matrices $\Hbf_d$, also referred to as faces or slices, contain the spatial information of the boundary element matrices at selected complex frequencies. The vectors $\f_d$, called fibers, contain the corresponding selected matrix entries evaluated over all frequencies $s_\ell$, $\ell=1,...,N_Q$. Hence, while the actual three-dimensional array $C$ would contain one face $\Hbf_d$ for each of the $N_Q$ frequencies, the 3D-ACA aims at representing this array by only $r$ face-fiber pairs, with $r \ll N_Q$. 
\begin{algorithm}
\caption{Pseudo code of 3D-ACA (slightly modified from \cite{seibel22})}\label{al:1}
\begin{algorithmic}[1]
\Function{$3$D ACA}{$\mathit{ENTRY}, \varepsilon$}
$C^{(0)} = 0, k_1=0 \textup{ and } d=0$
    \While{$\left\| \Hbf_{d} \right\|_F\left\| \f_{d} \right\|_2 > \varepsilon\left\| C^{(d)} \right\|_F $}
    	\State $d = d + 1$
        \State $H_{d}[i,j] = \mathit{ENTRY}(i,j,k_{d}) - C^{(d-1)}[i,j,k_{d}], \quad  i,j=1,\ldots,\frac{M}{3}, \quad k_{d}\in \{1,\ldots,N_Q\} $
        \Statex  \Comment{$H_{d}[i,j] \in \mathbb{C}^{3\times3}$}
        \State $H_d[i_d, j_d] = \max_{i,j} \sigma_{\min}\big(H_d[i,j]\big)$
        \State $f_{d}[k] = \big(H_{d}[i_{d},j_{d}]\big)^{-1} \big( \mathit{ENTRY}(i_{d},j_{d},k) - C^{(d-1)}[i_{d},j_{d},k] \big), \quad k=1,\ldots,N_Q$ 
        \Statex \Comment{$f_{d}[k] \in \mathbb{C}^{3\times3}$}
        \State $C^{(d)} = C^{(d-1)} + \mathbf{H_{d}\otimes f_{d}}$
        \State $k_{d+1} = \arg \max_{k} \sigma_{\min}\big(f_d[k]\big)$
    \EndWhile
    \State $r = d - 1$ \Comment{Final rank, \ie necessary frequencies}
    \State \Return $C^r = \displaystyle \sum_{d = 1}^r \Hbf_d\otimes \f_{d}$
\EndFunction
\end{algorithmic}
\end{algorithm}
The construction is performed iteratively, as outlined in Algorithm \ref{al:1}. Essentially, the algorithm takes as input the relevant spatial matrix at a selected frequency $s_\ell$, corresponding either to the collocation or to the Galerkin formulation, denoted by $\mathit{ENTRY}$, together with a prescribed precision $\varepsilon$. 

Starting from $d=1$, $\Hbf_d$ is filled with the input from $\mathit{ENTRY}$, after subtracting the current approximation $C^{\kl{d-1}}$, which is zero in the first iteration. The algorithm then searches for a pivot element $H_{d}[i_d,j_d]$. At the corresponding position $\kl{i_d,j_d}$, the fiber $\f_d$ is computed over all $N_Q$ frequencies and normalized by the pivot element $H_{d}[i_d,j_d]$. The outer product of this face-fiber pair is then added to the current approximation $C^{\kl{d-1}}$, and a pivot position is subsequently searched for in the fiber in an analogous way. At the corresponding position $k_{d+1}$, the next face is computed at the associated frequency. The iteration continues as long as the following accuracy condition holds:
\begin{equation}\label{eq:precision}
\left\| \Hbf_{d} \right\|_F\left\| \f_{d} \right\|_2 > \varepsilon\left\| C^{(d)} \right\|_F,    
\end{equation}
where the Frobenius norm $\left\| C^{(d)} \right\|_F$ does not require explicit assembly of the full array and can be computed recursively as \cite{seibel22}
\begin{equation}\label{eq:FrobNorm}
\left\| C^{(d)} \right\|_F^2 = \sum_{d,d'=1}^r\kl{\sum_{i,j}C_d[i,j]\overline{C_{d'}[i,j]}}\kl{\sum_k f_d[k] \overline{f_{d'}[k]}}.    
\end{equation}
Due to the tensorial nature of the elastodynamic problem, each element in $\Hbf$ is a $3\times 3$ matrix-block defined as
\begin{equation}
 H[i,j] := H[3i-2:3i, 3j-2:3j], \quad i,j=1,\ldots,\frac{M}{3},   
\end{equation}
and, hence, the pivot element block is chosen slightly different than in \cite{seibel22}. In the scalar case, the pivot position is simply chosen where the modulus attains its maximum value. In the present matrix valued case, however, the regularity of the pivot matrix must be guaranteed. Rjasanow and Weggler \cite{rjasanow17} considered several possibilities for this purpose and chose the following criterion for the selection of the pivot position, which is also adopted in the present work. Let the singular value decomposition (SVD) of $H_d[i,j]$ be
\begin{equation}\label{eq:SVD}
H_d[i,j] = U_d[i,j]\Sigma_d[i,j]V^*_d[i,j]    
\end{equation}
with 
\begin{equation*}
\Sigma_d[i,j] = \textup{diag}\kl{\sigma_1^{d}[i,j], \sigma_2^{d}[i,j], \sigma_3^{d}[i,j]}  
\end{equation*}
and
\begin{equation*}
\sigma_1^{d}[i,j] \geq \sigma_2^{d}[i,j] \geq \sigma_3^{d}[i,j] .
\end{equation*}
Then, the criterion for the choice of the pivot block $H_d[i_d,j_d]$ is given by
\begin{equation}
H_d[i_d,j_d] = \max_{i,j} \sigma_3\kl{H_d[i,j]}.
\end{equation}
Similarly, each fiber entry is a $3\times 3$ matrix-block,
\begin{equation*}
    f[k] \in \mathbb{C}^{3\times 3}.
\end{equation*}
Therefore, the pivot position in the fiber is selected in the same way as before, \ie the SVD is computed for each fiber entry, and the next frequency index $k_{d+1}$ is chosen as 
\begin{equation}
     k_{d+1} = \arg \max_{k} \sigma_{3}\big(f_d[k]\big).
\end{equation}
In the case that $\sigma_3\kl{H_d[i,j]} = 0$ for all $i,j$, or $\sigma_3\kl{f_d[k]} = 0$ for all $k$, no regular pivot matrix exists. In such a case, instead of using the so called Moore--Penrose pseudoinverse matrix as in \cite{rjasanow17}, the algorithm is skipped for that particular array $C$, which is then computed in dense form.  

The spatial discretization matrix $\Hbf$ is further approximated using either the $\mc{H}$-matrix based ACA or the Chebyshev interpolation based fast multipole method (FMM). 
\subsection{$\mc{H}$-matrix based ACA}
The spatial matrix approximation with $\mc{H}$-matrix based ACA is briefly outlined here; see, \eg \cite{bebendorf08} for more details. As an example, consider the Laplace transformed Galerkin single layer matrix $\hwmat{V}_D$. The Laplace transformed boundary element matrix $\hwmat{V}_D$ from \eqref{eq:discrGal} is given by 
\begin{equation}\label{eq:Vlap}
   \hwmat{V}_D[i,j]\kl{s_\ell} = \int_{supp\kl{\psi_i}}\psi_i\kl{\x}\int_{supp\kl{\psi_j}}\U\kl{\y -\x, s_\ell}\psi_j\kl{\y}\dsy\dsx = \Hbf[i,j]\kl{s_\ell}, \quad i\in I, j\in J,
\end{equation}
where $I,J = \{1,\ldots , M \}$. The matrix indices $i$ and $j$ are associated with the corresponding index sets $I$ and $J$, on the basis of which the matrix is decomposed into subblocks. Suitable subblocks are then approximated by low-rank representations. The partitioning of the matrix is performed by a recursive subdivision of the geometry. To this end, cluster trees $\mathscr{T}_I$ and $\mathscr{T}_J$ are constructed for $I$ and $J$ respectively.

Starting from the root level of $I$, the cluster $Cl_0\kl{I}$ is divided into two son clusters $Cl_1\kl{I}$ and $Cl_2\kl{I}$ by means of principal component analyses (PCA), such that
\begin{equation}\label{eq:cluster}
    Cl_1\kl{I} \cup Cl_2\kl{I} = Cl_0\kl{I}, \quad Cl_1\kl{I}\cap Cl_2\kl{I} = \emptyset.
\end{equation}
This subdivision is continued recursively until the cardinality of a cluster satisfies $\# Cl_x \leq b_{\textup{min}}$, at which point the cluster becomes a leaf cluster. Once the cluster trees $\mathscr{T}_I$ and $\mathscr{T}_J$ have been constructed in this manner, the block cluster tree $\mathscr{T}_{I\times J}$ is built from pairs of clusters $Cl_x\kl{I}\in \mathscr{T}_I$ and $Cl_y\kl{J}\in \mathscr{T}_J$, starting from the root block 
\begin{equation*}
Cl_0\kl{I} \times Cl_0\kl{J} = I\times J.
\end{equation*}
In this way, a hierarchical partition of the full matrix $\hwmat{V}_D\kl{s_\ell}$ into subblocks of the form $\kl{Cl_x\kl{I} \times Cl_y\kl{J}}$ is obtained, denoted by $\Hbf_{I\times J}$.

Since $\U\kl{\x -\y, s_\ell}$ in \eqref{eq:Vlap} is non-local and becomes singular as $\x \rightarrow \y$, not all subblocks can be approximated in low-rank. Therefore, all subblocks $\Hbf_{I\times J}$ are classified as near-field or far-field blocks according to the admissibility criterion
\begin{equation}\label{eq:adm}
    \textup{min} \{ \textup{diam}\kl{Cl_x\kl{I}},\textup{diam}\kl{Cl_y\kl{J}}  \} \leq \eta \textup{dist}\kl{Cl_x\kl{I},Cl_y\kl{J}},
\end{equation}
with a prescribed $\eta \in \mathbb{R}^+$. The subblocks satisfying this criterion are classified as far-field blocks, whereas the remaining blocks are classified as near-field blocks. The near-field matrix blocks are computed by standard BEM and stored densely, without approximation, while the far-field blocks are approximated by adaptive cross approximation (ACA). 

The matrix-valued ACA algorithm for the far-field subblocks is conceptually similar to the 3D-ACA algorithm in \ref{al:1}, but without the additional third dimension. A detailed explanation is therefore omitted here, and the reader is referred to \cite{haider19} for further details. Similar to the 3D-ACA, the low-rank representation of the spatial matrix subblock $\Hbf_{I\times J}$ is given by $\Hbf_{I\times J} \approx \Hbf^r_{I \times J} = \mathbf{U}\mathbf{V}^H$, where $\mathbf{U}$ and $\mathbf{V}$ are the rows and column matrices, respectively, analogous to the face-fiber pairs in 3D-ACA. Furthermore, the matrices $\mathbf{U}$ and $\mathbf{V}$ are subjected to a QR-decomposition, after which an SVD approximation is applied to the inner R-matrices of the decomposed factors. The complete low-rank approximation can be written as
\begin{equation}\label{eq:ACA}
\begin{aligned}
    \Hbf_{I\times J} \approx \Hbf^r_{I \times J} &= \mathbf{U}\mathbf{V}^H = \mathbf{Q}_{\mathbf{U}}\kl{\mathbf{R}_{\mathbf{U}}\mathbf{R}_{\mathbf{V}}^H}\mathbf{Q}_{\mathbf{V}}^H \\
    &= \mathbf{Q}_{\mathbf{U}} \kl{\check{\mathbf{U}}\check{\Sigma}\check{\mathbf{V}}^H} \mathbf{Q}_{\mathbf{V}}^H = \overline{\mathbf{U}}\check{\Sigma}\overline{\mathbf{V}}^H.
\end{aligned}
\end{equation}
However, the identification of the pivot matrix block in the approximated matrix subblock $ \Hbf^r_{I \times J} $ is not straightforward. Therefore, the corresponding matrix subblock is reconstructed from its low-rank approximation, the pivot matrix is determined, and the dense block is discarded again. Fortunately, these additional matrix operations can be implemented efficiently and have only a minor impact on the overall computation time. In the remainder of the paper, this ACA combined with the 3D-ACA approach is referred to simply as the "$\mc{H}$-matrix ACA" approach. 

\subsection{Fast multipole method (FMM)}
In the fast multipole method (FMM), the face matrix $\Hbf$ is never fully assembled. Instead, the matrix-vector multiplication is accelerated. In this sense, it is conceptually different from ACA and, therefore, requires a slightly different approach when combined with the 3D-ACA. The Chebyshev-interpolation based FMM \cite{fong09} is employed here as a black box with only minimal modifications. Standard geometric clustering techniques based on collocation points (or basis functions in Galerkin) are applied. Since the integration occurs over the supports of the basis functions, the bounding box of each cluster is extended so as to fully enclose these supports. A detailed description of the standard geometric clustering procedure can be found in the book by Liu \cite{liu09}. 

The core idea of the FMM is to compute a low-rank approximation of the integral kernels in the far-field. This kernel-dependent approximation shifts the dependence of the fundamental solution $\U$ on the distance function $|\x-\y|$ to two separate interpolating functions, solely dependent on either $\x$ or $\y$. Using an admissibility conditions similar to that of ACA in \eqref{eq:adm}, the cluster pairs are classified into far-field pairs, for which a low-rank approximation of the integral kernel is employed, and near-field pairs, for which the kernel is computed densely. Here, the low-rank approximation and the standard FMM operators are only briefly introduced. A detailed description of the FMM operators, as well as the interpolation and anterpolation techniques of the multilevel scheme, can be found in \cite{greengard87, fong09}. 

The low-rank approximation of the fundamental solution $\U\kl{\x,\y,s_\ell}$ is given by
\begin{equation}\label{eq:FMMU}
    \U\kl{\x,\y,s_\ell} \approx \sum_{\mathbf{n}}S_p\kl{\x,\bar{\x}_{\mathbf{n}}} \sum_{\mathbf{m}}\U\kl{\bar{\mathbf{x}}_{\mathbf{n}},\bar{\mathbf{y}}_{\mathbf{m}},s_\ell} S_p\kl{\bar{\y}_{\mathbf{m}},\y},
\end{equation}
where $\bar{\mathbf{x}}_{\mathbf{n}}$ and $\bar{\mathbf{y}}_{\mathbf{m}}$ are Chebyshev interpolation points. The interpolating function $S_p\kl{\x,\bar{\x}_{\mathbf{n}}}$ is given by
\begin{equation}\label{eq:int}
S_p\kl{\x,\bar{\x}_{\mathbf{n}}} := \Pi_{i=1}^3 S_p\kl{x_i,\bar{x}_{n_i}},    
\end{equation}
which is the three-dimensional extension of the standard one-dimensional Chebyshev interpolation function $S_p\kl{x_i,\bar{x}_{n_i}}$, given by
\begin{equation}\label{eq:int1}
S_p\kl{x_i,\bar{x}_{n_i}} = \frac{1}{p}+ \frac{2}{p}\sum_{k+1}^{p-1}T_k\kl{x}T_k\kl{\bar{x}_n} \quad \forall x \in [-1,1],
\end{equation}
where $T_k\kl{x}$ is the first-kind Chebyshev polynomial of order k and $\bar{x}_n$ are the corresponding roots. Inserting the low-rank approximation \eqref{eq:FMMU} into \eqref{eq:Vlap}, the low-rank approximation of the Laplace transformed Galerkin single layer matrix $\hwmat{V}_D[i,j]\kl{s_\ell}$ can be written as 
\begin{equation}\label{eq:FMMV}
\begin{aligned}
\hwmat{V}_D[i,j]\kl{s_\ell} &= \underbrace{\sum_{\mathbf{n}} \int_{supp\kl{\psi_i}}S_p\kl{\x,\bar{\x}_{\mathbf{n}}}\psi_i\kl{\x}}_{\textup{L2P-operator}} \; \; \underbrace{\sum_{\mathbf{m}} \U\kl{\y -\x, s_\ell}}_{\textup{M2L-operator}} \\ &\underbrace{\int_{supp\kl{\psi_j}}S_p\kl{\y,\bar{\y}_{\mathbf{m}}}\psi_j\kl{\y}}_{\textup{P2M-operator}}\dsy\dsx.
\end{aligned}
\end{equation}

The acronyms shown below the individual terms denote the standard abbreviations used in the multilevel FMM scheme. For the double layer and the hypersingular operators, the low-rank approximation is similar except for the presence of the traction operator $\mc{T}$. In that case, the traction operator is shifted onto the interpolating functions, which results in a tensor-valued interpolation operator, while the M2L-operator remains unchanged. As can be seen from \eqref{eq:FMMV}, only the M2L-operator depends on the frequency $s_\ell$ and is therefore considered in the 3D-ACA algorithm, making the identification of the pivot matrix block fairly straightforward. This FMM based 3D-ACA approach will be referred to as "FMM-ACA" in order to distinguish it clearly from the ACA based 3D-ACA approach.

Using either of the above low-rank approximations with the 3D-ACA algorithm, the multiplication in \eqref{eq:gCQ} is modified to 
\begin{equation}\label{eq:3DacaMult}
\begin{aligned}
\sum_{\ell=1}^{N_Q}\hmat{V}_{D}\kl{s_\ell}\mat{W}_\ell^{\Delta t_n} &= \sum_{\ell=1}^{N_Q}\sum_{k=1}^{r}\Hbf_k\kl{\hmat{V}_{D}} \otimes f_k[\ell] \kl{\hmat{V}_{D}}\mat{W}_\ell^{\Delta t_n}    \\
&= \sum_{k=1}^{r}\Hbf_k\kl{\hmat{V}_{D}}\otimes \sum_{\ell=1}^{N_Q}f_k[\ell] \kl{\hmat{V}_{D}}\mat{W}_\ell^{\Delta t_n} .
\end{aligned}
\end{equation}
The complexity is thereby reduced from $\bigO\kl{M^2N_Q}$ in \eqref{eq:gCQ} to $\bigO\kl{r\kl{M^2+N_Q}}$ in \eqref{eq:3DacaMult}. The main factor contributing to this reduction is that the outer sum, which corresponds to a matrix vector multiplication of size $M$, is performed only $r$ times, instead of $N_Q$ as in the dense case. Hence, by selecting the frequencies greedily and approximating the resulting three-dimensional data array, the 3D-ACA algorithm enables a substantial reduction in computational complexity.


\section{Numerical examples}
\label{sec5}

The two proposed approaches, $\mc{H}$-matrix-ACA and FMM-ACA, are first tested on a unit cube subjected to a smooth pulse. The aim is to demonstrate that the low-rank approximations in both the spatial and the frequency dimensions do not affect the results. To this end, a convergence study is carried out for both a pure Dirichlet problem using collocation BEM and a mixed boundary value problem using Galerkin BEM. The corresponding storage compression rates of the two approaches are then compared, with particular emphasis on the number of frequencies selected by the Algorithm \ref{al:1} and on the regions of the complex plane in which most of these frequencies are located. As a second example, the vibrations of an electric machine are investigated by means of two different test cases. Again, the two proposed approaches are compared in terms of their compression rates and computation times. In all examples, the geometry is discretized using linear triangular elements. For the time stepping scheme, a 2-stage Radau IIA method is employed, and the resulting systems of equations are solved using BiCGstab without a preconditioner. 

\subsection{Unit cube subjected to a smooth pulse}\label{sec:cube}

\begin{figure}[!h]
    \centering
    \begin{subfigure}[b]{0.45\textwidth}
        \includegraphics[scale = 0.13]{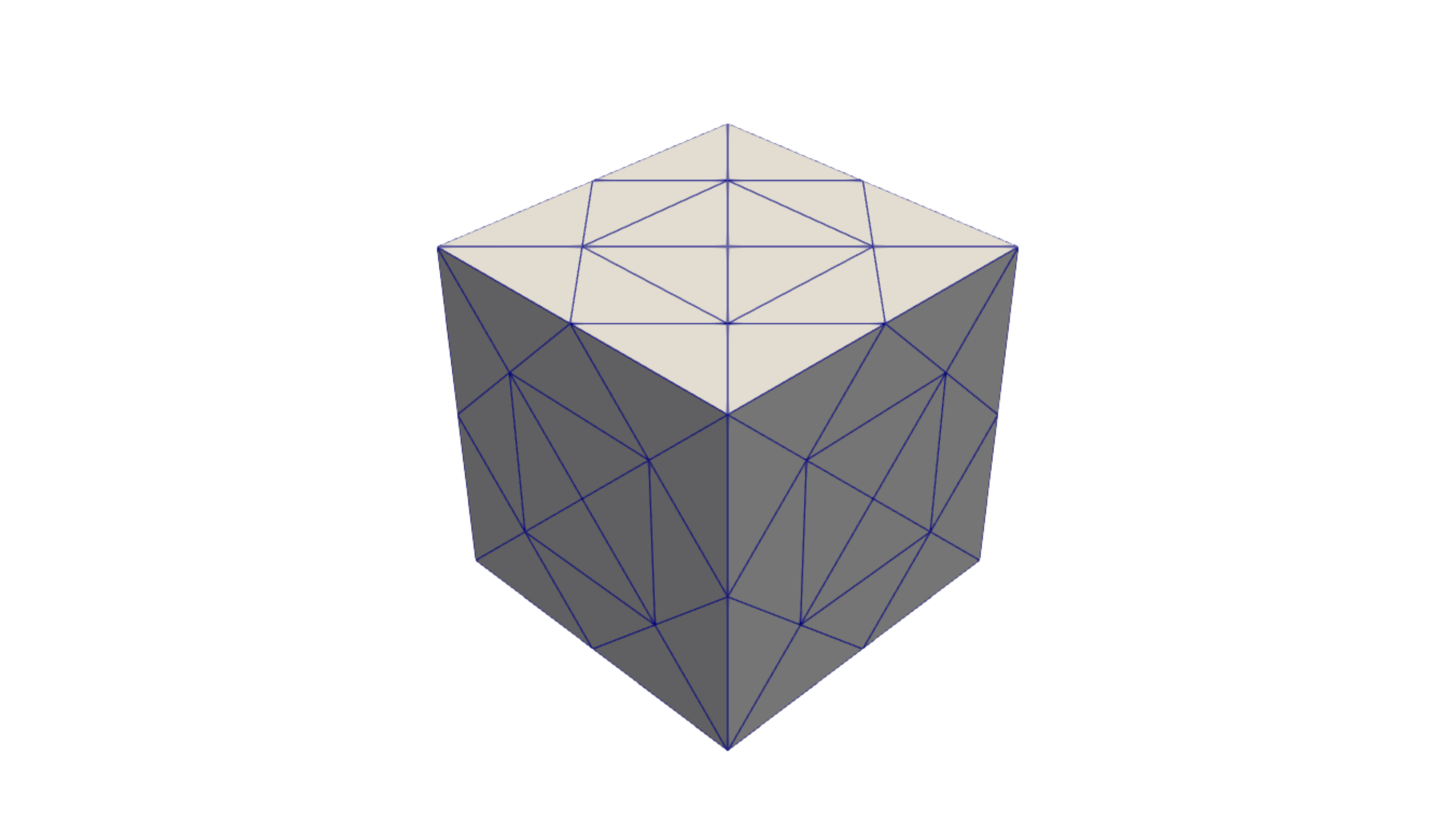}
        \caption{lvl 1: 96 elements}
        \label{fig:cube1}
    \end{subfigure}
    \begin{subfigure}[b]{0.45\textwidth}
        \includegraphics[scale = 0.13]{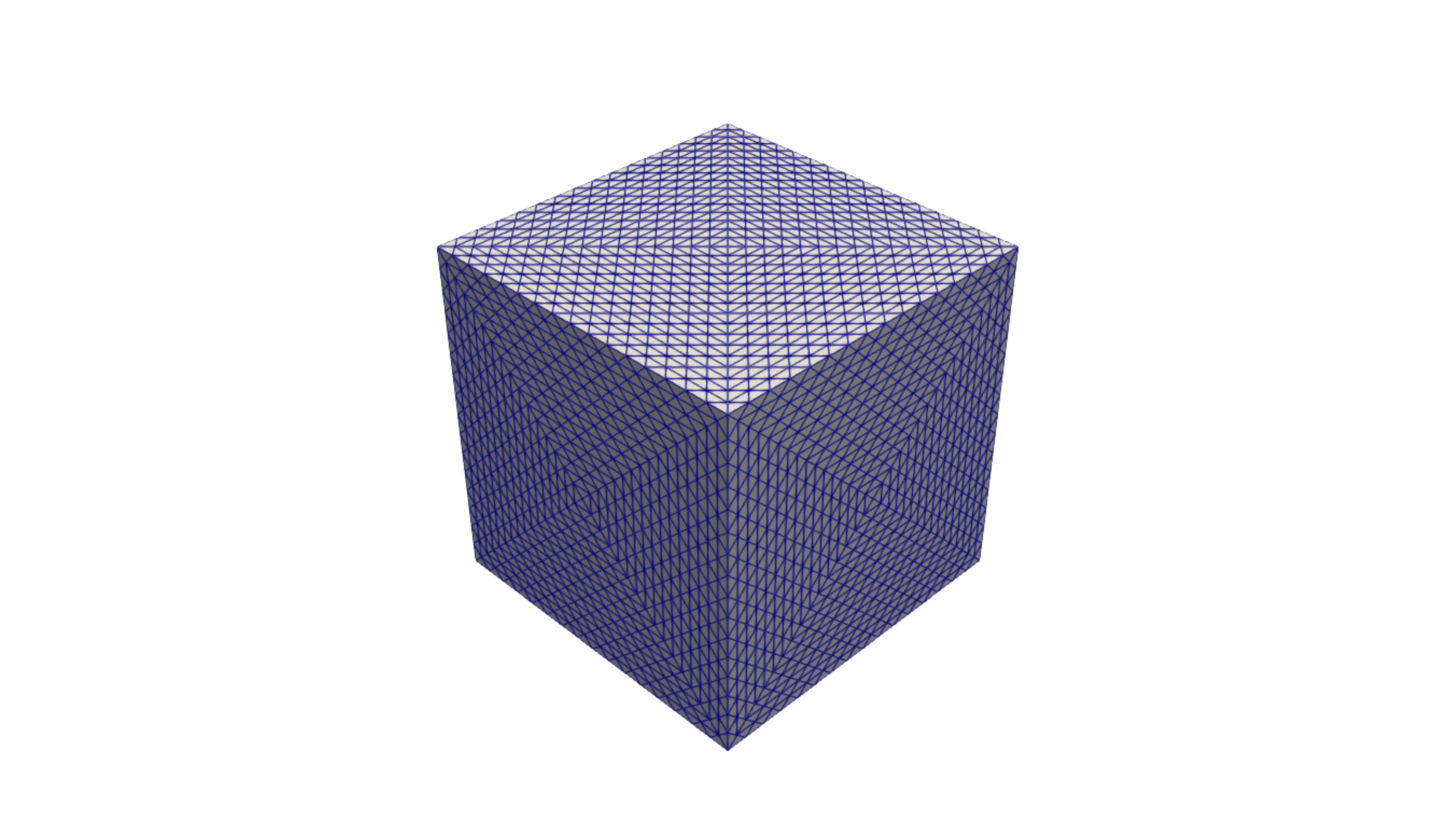}
        \caption{lvl 4: 6144 elements}
        \label{fig:cube2}
    \end{subfigure}
    \caption{Unit cube}
    \label{fig:cube}
\end{figure}

A unit cube $[-0.5,0.5]^3$ centred at the origin is considered, as shown in Fig. \ref{fig:cube}. 
A convergence study is done for the Dirichlet problem solved by collocation BEM and a mixed boundary value problem solved by the Galerkin BEM formulation. The initial unrefined mesh (Fig. \ref{fig:cube1}) consists of 96 elements. After four levels of refinement, the finest mesh (Fig. \ref{fig:cube2}) consists of 6144 elements. The refinement levels and the corresponding mesh data are listed in Table \ref{tab1}, where $h$ denotes the mesh size and $\Delta t$ the time-step size, while $N$ is the number of time steps and $N_Q$ is the number of complex frequencies used in the gCQ.   
\begin{table}[]
    \centering
    \begin{tabular}{l|rrrrrr}
    \toprule
       refinement  & nodes & elements & $h$ & $\Delta t$ & $N$ & $N_Q$ \\
    \midrule
        $1$ & $50$  & $96$ & $\SI{0.5}{\metre}$ & $\SI{0.3}{\second}$ & $10$ & $27$ \\
        $2$ & $194$ & $384$ & $\SI{0.25}{\metre}$ & $\SI{0.15}{\second}$ & $20$ & $90$ \\
        $3$ & $770$ & $1536$ & $\SI{0.125}{\metre}$ & $\SI{0.075}{\second}$ & $40$ & $272$ \\
        $4$ & $3074$ & $6144$ & $\SI{0.0625}{\metre}$ & $\SI{0.0375}{\second}$ & $80$ & $768$ \\ 
    \bottomrule
    \end{tabular}
    \caption{Mesh details for different refinement levels}
    \label{tab1}
\end{table}
For this example, the problems are investigated using the material parameters
\begin{equation}
    c_1=\SI{1}{\metre\per\second}, \quad c_2=\sqrt{\frac{1}{2}}\,\si{\metre\per\second}, \quad \rho=\SI{1}{\kilogram\per\cubic\metre},
\end{equation}
where $c_1$ and $c_2$ are the compression wave and shear wave speed respectively. 

The tolerance $\varepsilon$ for the stopping criterion in Algorithm \ref{al:1}, as well as for the ACA tolerance for the faces, $\varepsilon_{\textup{ACA}}$, and the corresponding FMM levels and the interpolation orders, are listed in Table \ref{tab2}. 
\begin{table}[]
    \centering
    \begin{tabular}{l|llll}
    \toprule
        refinement & $\varepsilon$ & $\varepsilon_\textup{ACA}$ & $\textup{F}_{lvl}$ & $p$ \\
    \midrule
        $1$ & $10^{-2}$ & $10^{-4}$ & $1$ & $2$ \\
        $2$ & $10^{-3}$ & $10^{-5}$ & $2$ & $3$ \\
        $3$ & $10^{-4}$ & $10^{-6}$ & $2$ & $4$ \\
        $4$ & $10^{-5}$ & $10^{-7}$ & $3$ & $5$ \\
    \bottomrule
    \end{tabular}
    \caption{Parameters for the $\mc{H}$-matrix ACA and FMM-ACA}
    \label{tab2}
\end{table}
The tolerance $\varepsilon$ for the 3D-ACA is selected as $\varepsilon = 100\varepsilon_{\textup{ACA}}$, where both tolerances are reduced by a factor of $10$ for each refinement level. The quantities $\textup{F}_{lvl}$ and $p$ denote the number of levels in the multilevel FMM scheme and the interpolation order, respectively. 

\begin{figure}[!h]
    \centering
    \begin{subfigure}[b]{0.4\textwidth}
        \includegraphics[scale = 1]{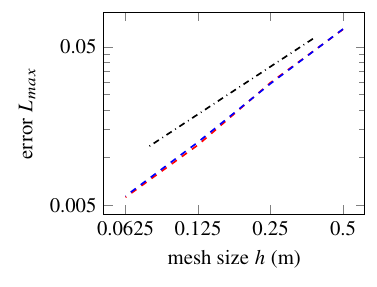}
        \caption{Dirichlet problem}
        \label{fig:conv1}
    \end{subfigure}
    \begin{subfigure}[b]{0.49\textwidth}
        \includegraphics[scale = 1]{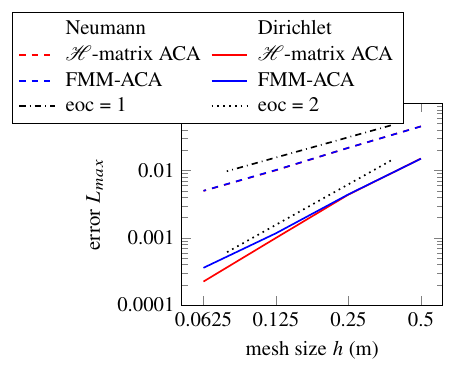}
        \caption{Mixed problem}
        \label{fig:conv2}
    \end{subfigure}
    \caption{Unit cube: $L_{max}$-error versus refinement in space and time}
    \label{fig:conv}
\end{figure}

For the prescription of the boundary conditions, an analytic full-space displacement field $\ubf\kl{\x,t}$ and traction field $\tbf\kl{\x,t}$ are constructed according to \cite{eringen74}, with the expressions listed in appendix \ref{apndxA}. This construction requires the definition of a source point $\mathbf{p}$ and a direction vector $\mathbf{d}$, chosen here as $\mathbf{p}= [1.5,1.2,1.2]^\top$ and $\mathbf{d} = [1,1,1]^\top$. In the Dirichlet case, $\ubf\kl{\x,t}$ is applied on the entire boundary. In the mixed problem, one half of the cube is prescribed with $\ubf\kl{\x,t}$, while the other half is prescribed with $\tbf\kl{\x,t}$. The error analysis is based on a uniform refinement in both $h$ and $\Delta t$, carried out simultaneously in space and time. The combined space-time error is measured by the $L_2$-error in space, and the maximum value over all time steps is taken. It is computed as 
\begin{equation}\label{eq:L2}
   L_{max}(\ubf_h) = \max_{1 \leq n \leq N} \left\| \ubf \kl{\frac{t_n + t_{n+1}}{2} } - \ubf_h \kl{ \frac{t_n + t_{n+1}}{2} } \right\|_{L_2},
\end{equation}
where $\ubf_h$ and $\ubf$ denote the computed quantity and the corresponding analytical solution, respectively. The order of convergence ($\textup{eoc}$) is then given by
\begin{equation}\label{eq:eoc}
     \textup{eoc} =  \log_2 \kl{ \frac {L_{max}^k} {L_{max}^{k+1} }},
\end{equation}
where $k$ denotes the corresponding refinement level. Although the underlying 2-stage Runge--Kutta yields third-order convergence in time, the combined error in space and time is dominated by the lower order spatial shape functions. Therefore, the Neumann data is expected to exhibit a first-order convergence, whereas the Dirichlet data is expected to exhibit quadratic order of convergence. It is worth mentioning the work of Kramer \etal \cite{kramer26}, where higher-order convergence in both space and time was achieved by combining CQ with isogeometric analysis. The error plots for both the Dirichlet and the mixed problems are shown in Figure \ref{fig:conv}.
The results show a very good order of convergence for both the Dirichlet and the mixed problems, and for both proposed approaches. In particular, the error curves overlap almost completely for the coarser meshes, while the FMM-ACA approach performs only slightly worse than its counterpart on the finer levels. However, this difference is mainly due to the chosen FMM parameters, which could be further improved, for example, by increasing the interpolation order $p$. 

\begin{figure}[!h]
    \centering
    \begin{subfigure}[b]{0.4\textwidth}
        \includegraphics[scale = 1]{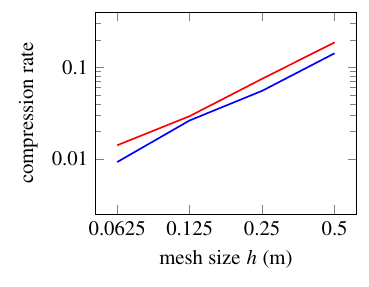}
        \caption{Dirichlet problem}
        \label{fig:comp1}
    \end{subfigure}
    \begin{subfigure}[b]{0.49\textwidth}
        \includegraphics[scale = 1]{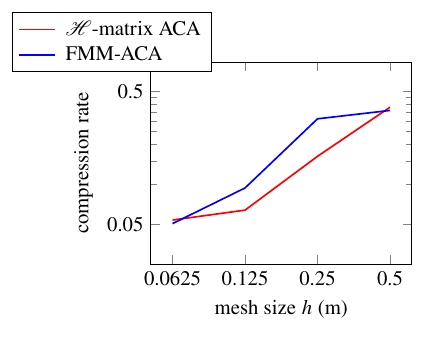}
        \caption{Mixed problem}
        \label{fig:comp2}
    \end{subfigure}
    \caption{Compression rates versus refinement for both problems}
    \label{fig:comp}
\end{figure}

\begin{figure}[!h]
    \centering

    \begin{subfigure}[b]{0.45\textwidth}
        \centering
        \includegraphics[scale=1]{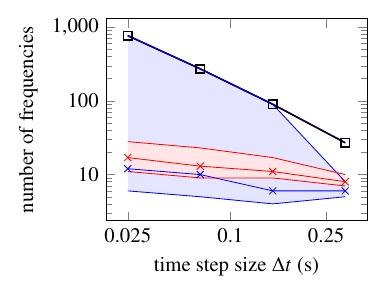}
        \caption{Dirichlet problem: Single-layer potential}
        \label{freq:dirislp}
    \end{subfigure}
    \hfill
    \begin{subfigure}[b]{0.45\textwidth}
        \centering
        \includegraphics[scale=1]{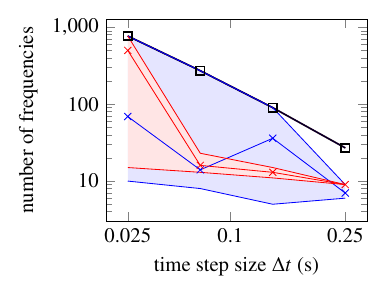}
        \caption{Mixed problem: Single-layer potential}
        \label{freq:slp}
    \end{subfigure}

    \vspace{0.4cm}

    \begin{subfigure}[b]{0.46\textwidth}
        \centering
        \includegraphics[scale=1]{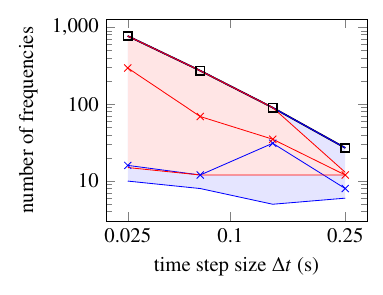}
        \caption{Mixed problem: Double-layer potential}
        \label{freq:dlp}
    \end{subfigure}
    \hfill
    \begin{subfigure}[b]{0.45\textwidth}
        \centering
        \includegraphics[scale=1]{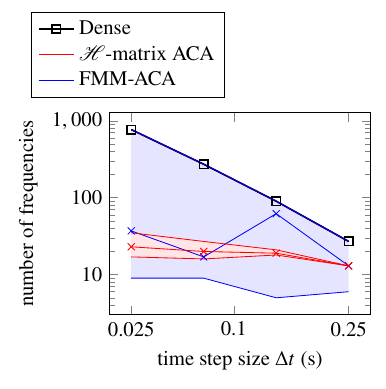}
        \caption{Mixed problem: Hypersingular operator}
        \label{freq:hso}
    \end{subfigure}

    \caption{Unit cube: used number of frequencies for the different BEM matrices for both problems}
    \label{fig:freq}
\end{figure}

More interesting are the compression rates achieved by the two approaches with respect to a dense BEM formulation. Figure \ref{fig:comp} shows the compression rates versus the refinement level for both problems. 
The compression rates shown here are defined as the memory used by the $\mc{H}$-matrix ACA or the FMM-ACA approach to the memory required by the corresponding fully dense computation for storing the matrices in \eqref{eq:discrCollo} and \eqref{eq:discrGal}. The results clearly show that high compression rates are achieved for both proposed approaches, independently of the problem under consideration and of the formulation employed. The FMM-ACA approach yields slightly better compression, especially at higher refinement levels, which is consistent with the observations in the acoustic case \cite{schanz26}. The main contribution to these compression rates comes from the 3D-ACA algorithm through its low-rank approximation of the three-dimensional data array. Consequently, the rank reduction with respect to the frequencies is the next aspect of interest. Since geometric clustering is employed in both the $\mc{H}$-matrix ACA and FMM-ACA approaches, the 3D-ACA algorithm is applied separately to each individual cluster block. As a result, the rank varies from block to block within the boundary element matrix. Figure \ref{fig:freq} therefore reports the mean number of retained frequencies together with a min--max envelope. The shaded band represents the variation of the retained frequency ranks across all blocks. In Figures \ref{freq:dirislp} and \ref{freq:slp}, the number of used frequencies for the single-layer potential is shown for both problems. 
Both the $\mc{H}$-matrix ACA and the FMM-ACA approaches show very low minimum and mean numbers of used frequencies across all cluster blocks. However, some blocks, especially in the FMM-ACA approach, had to be treated as dense blocks due to the absence of a regular pivot matrix, and hence the maximum reaches the dense limit. Figures \ref{freq:dlp} and \ref{freq:hso} show the used frequencies in the double-layer potential and the hypersingular operator of the mixed problem, respectively. The double-layer potential stands out somewhat from the others, since for both approaches the maximum number of used frequencies reaches the dense limit, for the same reason as before. Nevertheless, the mean values indicate that a large number of blocks are still compressed, with the FMM-ACA performing slightly better. Figure \ref{freq:hso} appears to suggest that the $\mc{H}$-matrix ACA approach is more suitable overall, although the mean values remain very similar for both approaches. 

It should be noted, however, that the comparison between the two approaches is not entirely direct. In the FMM-ACA approach, only the M2L blocks are considered in the far-field blocks, and these correspond to the single-layer operator. Consequently, for the FMM-ACA approach, the data shown in Figure \ref{fig:freq} reflects only the near-field contribution, resulting in the comparison to be slightly biased. Overall, the mean number of frequencies selected by the 3D-ACA algorithm remains low, irrespective of the problem type, the approach employed, or the operator considered. This continues to be the major contributor to the compression of the three-dimensional data array. 

The complex frequencies at which the spatial matrices are evaluated are determined by the gCQ method and lie on a circle in the right half of the complex plane, whose radius and position depend on certain parameters. Since the same parameters as in \cite{schanz26} are used here, they are not repeated. These complex frequencies occur in conjugate pairs, so that in a dense computation only one half of them needs to be evaluated explicitly.
\begin{figure}[!h]
    \centering

    \begin{subfigure}[b]{0.45\textwidth}
        \centering
        \includegraphics[scale=1]{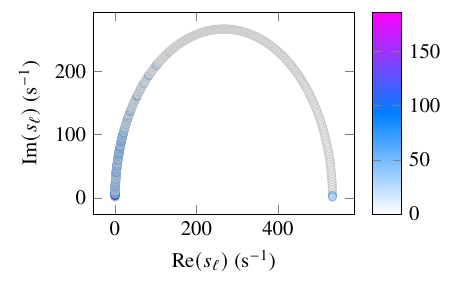}
        \caption{$\mc{H}$-matrix ACA SLP}
        \label{fig:ACAslp}
    \end{subfigure}
    \hfill
    \begin{subfigure}[b]{0.45\textwidth}
        \centering
        \includegraphics[scale=1]{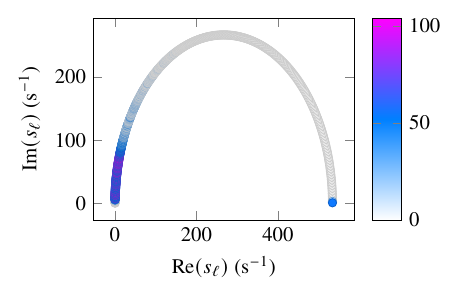}
        \caption{FMM-ACA SLP}
        \label{fig:FMMslp}
    \end{subfigure}

    \vspace{0.4cm}

    \begin{subfigure}[b]{0.45\textwidth}
        \centering
        \includegraphics[scale=1]{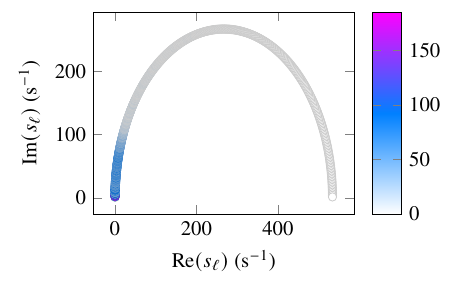}
        \caption{$\mc{H}$-matrix ACA DLP}
        \label{fig:ACAdlp}
    \end{subfigure}
    \hfill
    \begin{subfigure}[b]{0.45\textwidth}
        \centering
        \includegraphics[scale=1]{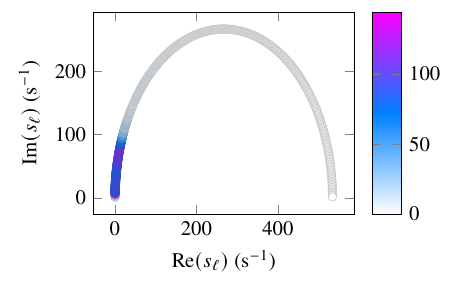}
        \caption{FMM-ACA DLP}
        \label{fig:FMMdlp}
    \end{subfigure}

        \vspace{0.4cm}

    \begin{subfigure}[b]{0.45\textwidth}
        \centering
        \includegraphics[scale=1]{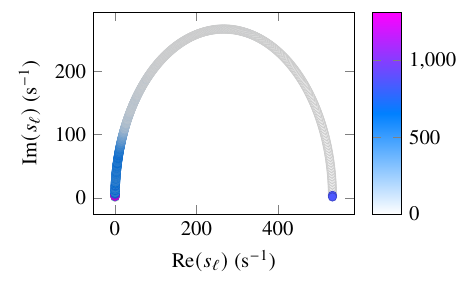}
        \caption{$\mc{H}$-matrix ACA HSO}
        \label{fig:ACAhso}
    \end{subfigure}
    \hfill
    \begin{subfigure}[b]{0.45\textwidth}
        \centering
        \includegraphics[scale=1]{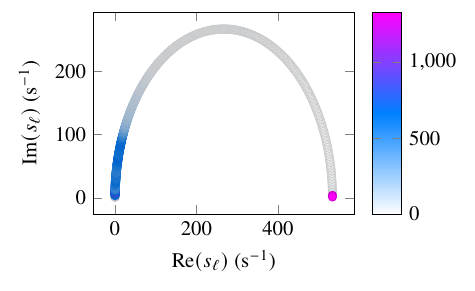}
        \caption{FMM-ACA HSO}
        \label{fig:FMMhso}
    \end{subfigure}

    \caption{Unit cube (refinement 4): used complex frequencies. The colour code shows the number of cluster blocks that uses the corresponding frequency}
    \label{fig:circle}
\end{figure}
Figure \ref{fig:circle} shows the frequencies selected by the 3D-ACA algorithm on the half-circle, with a colour coding indicating the number of cluster blocks that select the corresponding frequency. Note that dense cluster blocks are excluded from this plot. Moreover, for those cluster blocks that are compressed by the 3D-ACA, the first two frequencies are omitted from the colour plots, since they are selected by almost all blocks. From Figure \ref{fig:circle}, it is evident that Algorithm \ref{al:1} predominantly selects frequencies with small real parts $\textup{Re}\kl{s_\ell}$, irrespective of the approach or the operator considered. Interestingly, for the single-layer and the hypersingular operators, the algorithm also selects the last few frequencies on the half-circle, \ie frequencies with a large $\textup{Re}\kl{s_\ell}$. 

\begin{figure}[!h]
    \centering
    \begin{subfigure}[b]{0.45\textwidth}
        \includegraphics[scale = 1]{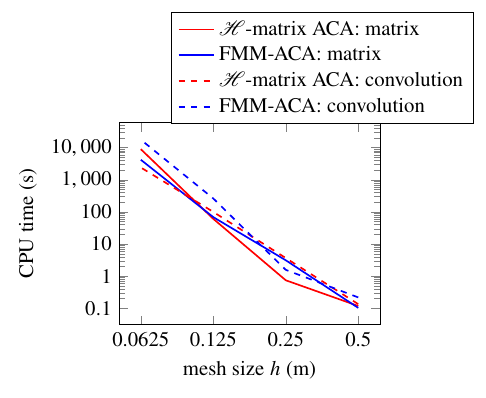}
        \caption{Dirichlet problem: collocation BEM}
        \label{fig:diri_time}
    \end{subfigure}
    \begin{subfigure}[b]{0.45\textwidth}
        \includegraphics[scale = 1]{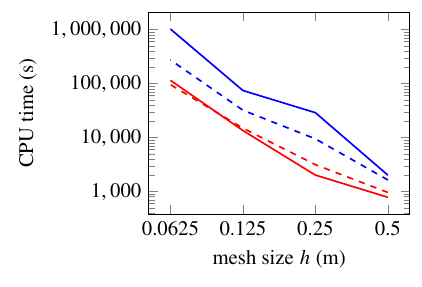}
        \caption{Mixed problem: Galerkin BEM}
        \label{fig:mixed_time}
    \end{subfigure}
    \caption{Unit cube: Computing time for assembly and solving}
    \label{fig:time}
\end{figure}

Finally, Figure \ref{fig:time} shows the computing times versus the refinement levels for both problems.
Essentially, two different contributions to the total computing time are compared. The first is the time required by the 3D-ACA algorithm to construct the data-sparse array, denoted as the 'matrix' time. The second is the convolution in \eqref{eq:3DacaMult}, denoted as the 'convolution' time. It should be noted that the mixed problem is solved using the Galerkin formulation, which leads to significantly larger CPU times than in the Dirichlet problem. For the Dirichlet problem, it is evident from Figure \ref{fig:diri_time} that the matrix assembly with FMM-ACA becomes much faster than with $\mc{H}$-matrix ACA at the finer refinement levels, for two reasons. At the coarser levels, the multilevel FMM scheme consists only of the near-field and, therefore, all the operator matrices has to be computed by the 3D-ACA, and this takes longer than in the $\mc{H}$-matrix ACA approach. At finer levels, however, the far-field computation involves only the M2L-operator and the fibers, and is therefore much faster. For the mixed problem shown in Figure \ref{fig:mixed_time}, on the other hand, the matrix assembly takes considerably longer with the FMM-ACA approach, in contrast to the collocation BEM. Although the FMM-ACA approach involves only the M2L-operator in the far-field, it was observed that the near-field hypersingular operator becomes the bottleneck, leading to much longer assembly times than in the $\mc{H}$-matrix ACA approach. As for the convolution times, both problems show that the FMM-ACA approach is considerably slower, irrespective of whether collocation or Galerkin BEM is employed. Overall, although the  FMM-ACA approach yields very similar errors and slightly better compression rates, it requires noticeably more total computing time than the $\mc{H}$-matrix ACA. It should be emphasized, however, that both the chosen parameters and the uniform cluster tree in the FMM plays a crucial role, and that the results might be improved by selecting a different set of parameters or by using a more balanced cluster tree.   

\subsection{Vibrations of an electric machine}
In this section, two different test cases for analysis of vibrations in an electric machine are presented. The main objective is to demonstrate the efficiency and robustness of the two proposed approaches for large-scale problems. In the first test case, a localized impact hammer test is modelled by applying a Gamma-shaped pulse on a very small area, and the resulting vibration response is analysed in the time-domain. In the second test case, a quadrupole excitation is introduced in the radial direction, in order to mimic the Maxwell stresses acting on the interior of the electric machine, and the resulting vibrations of the exterior casing are analysed. Both problems are solved by collocation BEM using the two proposed approaches, and a comparison is presented.

\begin{figure}[!h]
    \centering
    \begin{subfigure}[b]{0.6\textwidth}
        \includegraphics[scale = 0.3,
        trim={18cm 9cm 18cm 10cm},
        clip]{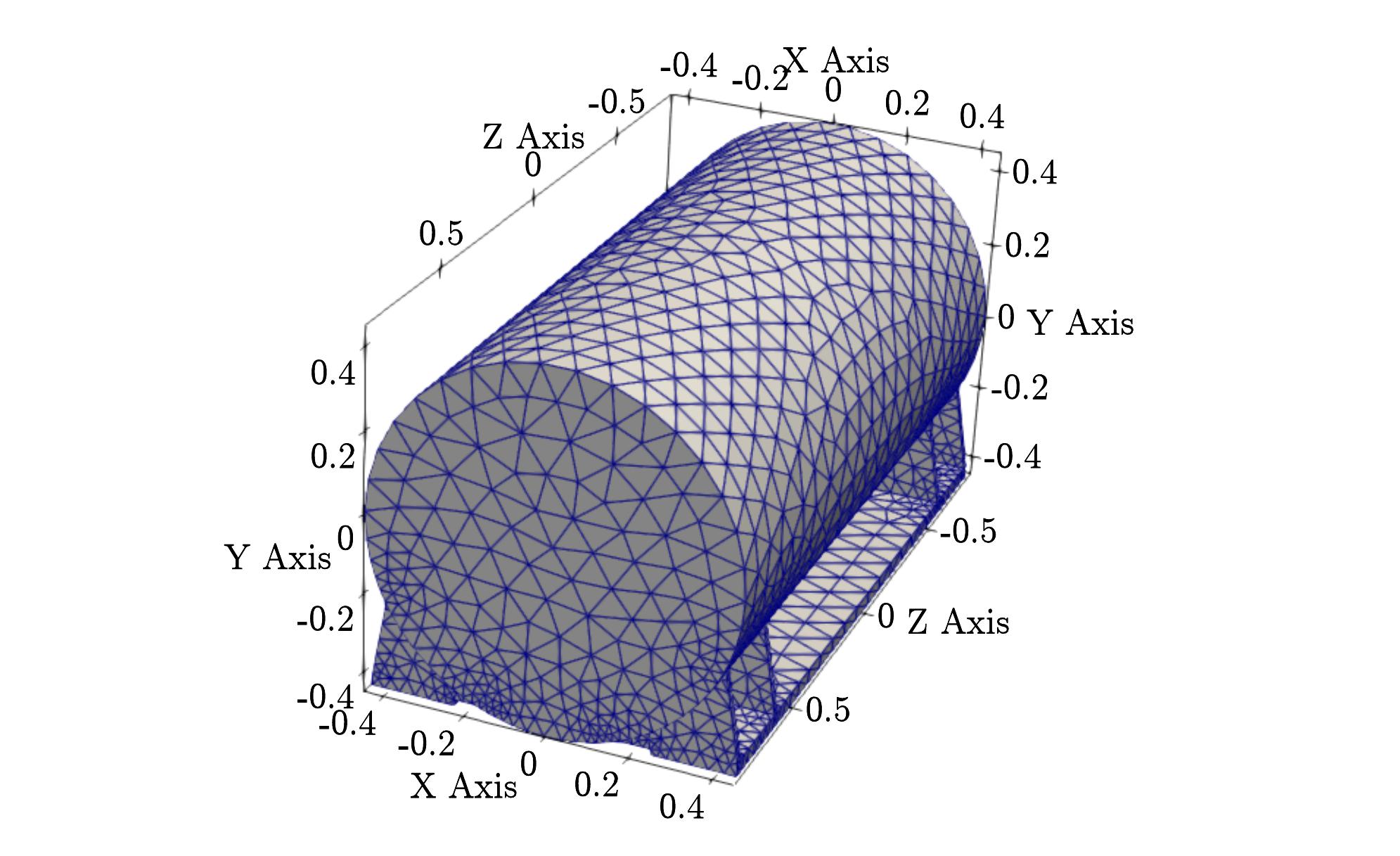}
        \caption{Electric machine}
        \label{fig:machineA}
    \end{subfigure}
    \begin{subfigure}[b]{0.3\textwidth}
        \begin{tabular}{ll}
            nodes &  $3363$\\
            elements & $6718$\\
            mean $h$ &  $\SI{0.06}{\metre}$ \\
            $\Delta t_1$, $\Delta t_2$ & $5\cdot 10^{-5}\,\si{\second}$, $10^{-4}\,\si{\second}$ \\
            $N_1$, $N_2$  & $120$, $350$ \\
            $\varepsilon_{\textup{ACA}}$ & $10^{-7}$ \\
            $\varepsilon$ & $10^{-5}$ \\
            $\textup{F}_{lvl}$, $p$ & $3$, $5$ \\
            & \\
        \end{tabular}
        \caption{Computing data}
    \label{fig:machine_data}
    \end{subfigure}
    \caption{Electric machine geometry and parameters}
    \label{fig:machine}
\end{figure}

The electric machine considered here is similar to the induction machine model studied in \cite{weilharter11}. The meshed geometry is shown in Figure \ref{fig:machine}, together with the corresponding spatial and temporal parameters. Here, subscripts $_1$ and $_2$ refer to the first and second test cases, respectively.
Since the quantity of interest is the vibration of the exterior casing, the interior stator and rotor assembly is neglected, and only the exterior casing is considered, having a thickness of $\SI{0.15}{\metre}$ and is assumed to be made of aluminium, with material parameters
\begin{equation}
    c_1 =  \SI{5994.26}{\metre\per\second}, \quad c_2 =  \SI{3204.07}{\metre\per\second}, \quad \nu = 0.3, \quad \rho = \SI{2660}{\kilogram\per\cubic\metre}.
\end{equation}

\subsubsection{Impact hammer test}
\begin{figure}
    \centering
    \includegraphics[scale=1]{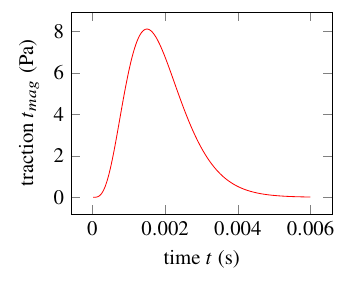}
    \caption{Input traction}
    \label{fig:traction}
\end{figure}

In this test case, the impact of a hammer on the electric machine is modelled by prescribing a smooth Gamma-shaped pulse as an excitation on a very small area of the machine casing. The pulse is given by
\begin{equation}\label{eq:smoothpulse}
    f\kl{t} = \kl{at}^m \textup{exp}\kl{-at},
\end{equation}
where the parameters $m$ and $a$ are chosen as $m=4$ and $a = \frac{4m}{N_1\Delta t}$, such that the wave can propagate across the machine completely. A mixed boundary value problem is considered, with the base of the machine prescribed with homogeneous Dirichlet conditions. A small area consisting of 6 boundary elements on the rim of the casing is prescribed with the traction boundary condition 
\begin{equation}\label{eq:machineBC1}
   \tbf\kl{\x,t} =  f\kl{t} \n\kl{\x},
\end{equation}
while the remaining part of the boundary is prescribed with zero traction. Figure \ref{fig:traction} shows the magnitude of the prescribed traction at one boundary element.

\begin{figure}[!h]
    \centering
    \begin{subfigure}[b]{0.49\textwidth}
        \includegraphics[scale = 0.2,
        trim={7cm 2cm 5cm 2cm},
        clip]{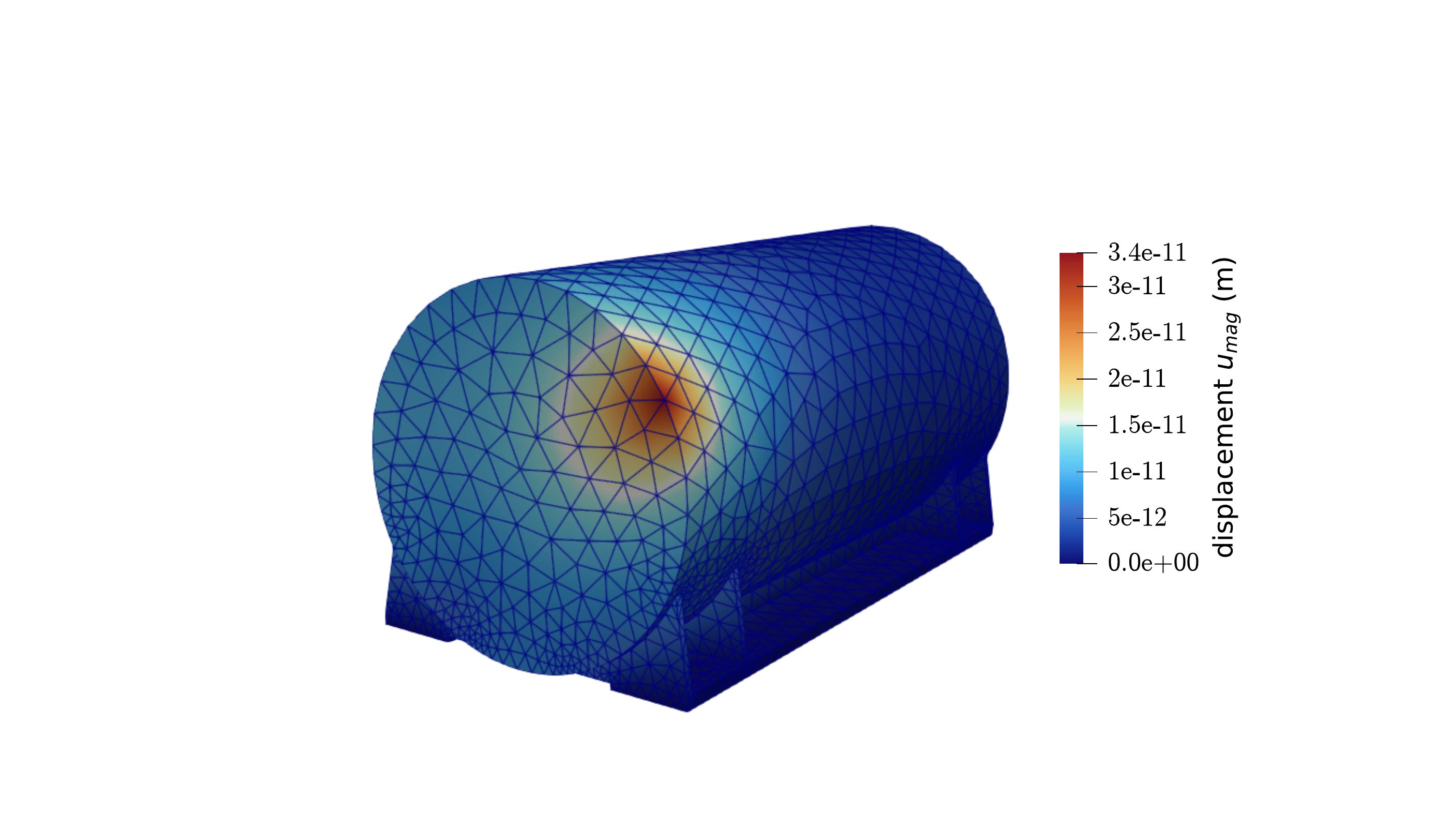}
        \caption{Time $t =  \SI{0.00100}{\second}$}
        \label{fig:hammer1}
    \end{subfigure}
    \begin{subfigure}[b]{0.45\textwidth}
        \includegraphics[scale = 0.2,
        trim={7cm 2cm 5cm 2cm},
        clip]{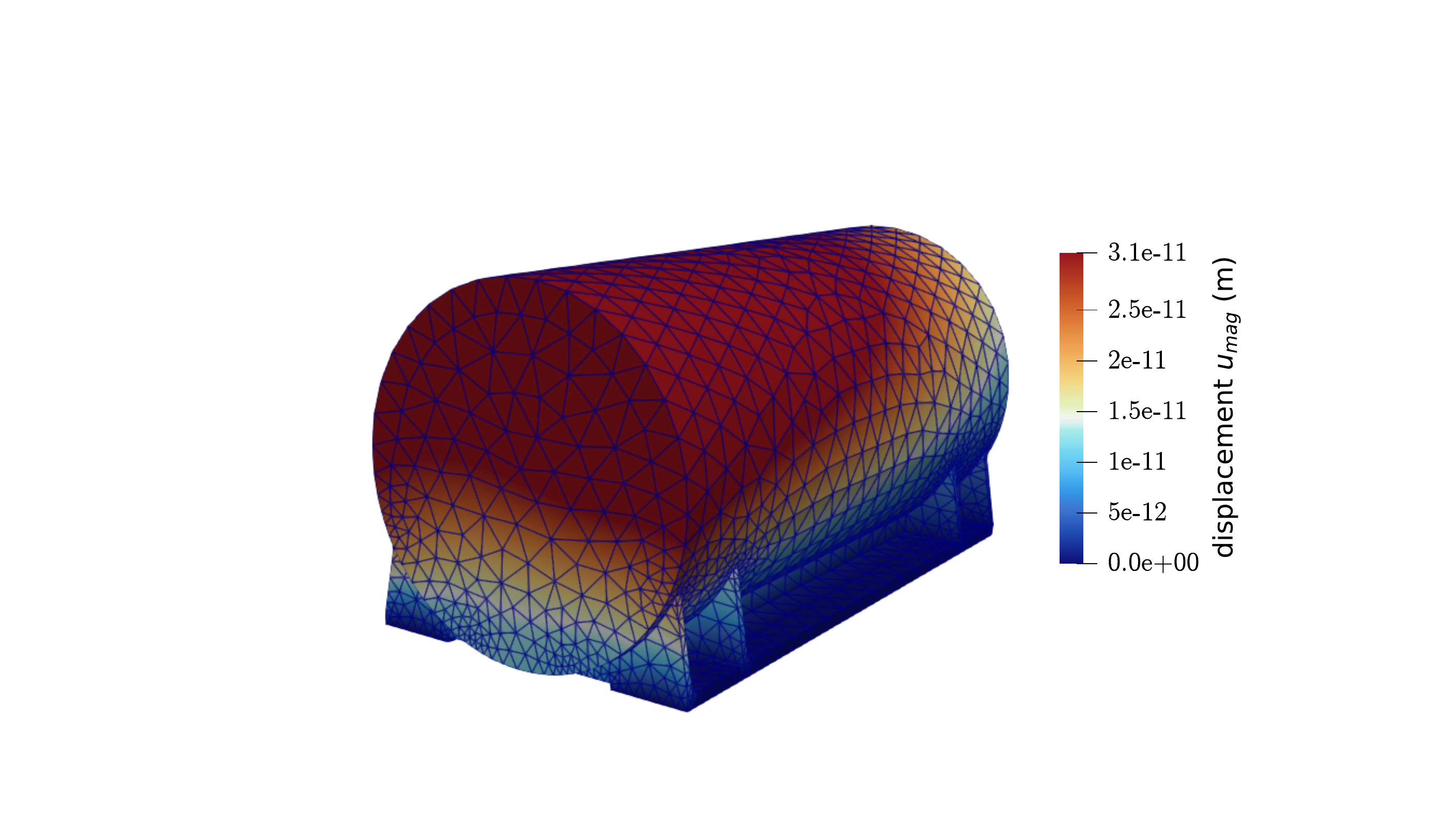}
        \caption{Time $t =  \SI{0.00202}{\second}$ }
        \label{fig:hammer2}
    \end{subfigure}
    \caption{Magnitude of displacement $u_{mag}$ (in m) on the surface of the machine casing at different times ($\mc{H}$-matrix ACA)}
    \label{fig:hammer}
\end{figure}

\begin{figure}[!h]
    \centering
    \begin{subfigure}[b]{0.45\textwidth}
        \includegraphics[scale = 1]{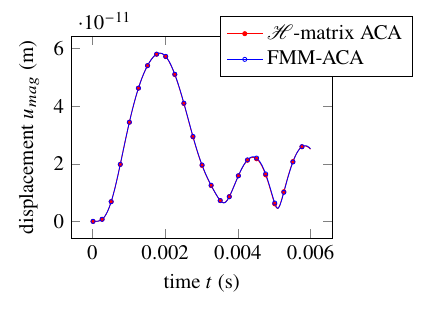}
        \caption{$u_{mag}$ at node A}
        \label{fig:disp_hammer1}
    \end{subfigure}
    \begin{subfigure}[b]{0.45\textwidth}
        \includegraphics[scale = 1]{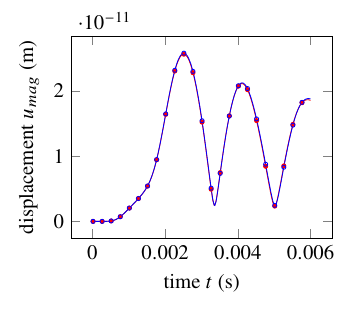}
        \caption{$u_{mag}$ at node B}
        \label{fig:disp_hammer2}
    \end{subfigure}
    \caption{Magnitude of displacement $u_{mag}$ over time at two distinct nodes}
    \label{fig:disp_hammer}
\end{figure}

First, the results obtained with the $\mc{H}$-matrix ACA approach are presented in Figure \ref{fig:hammer}.
The displacement magnitude $u_{mag}$ is shown on the surface of the machine casing at two different time instants. Figure \ref{fig:hammer1} shows the displacement at time $t=\SI{0.001}{\second}$, slightly before the peak magnitude is reached, whereas \ref{fig:hammer2} shows the displacement at time $t= \SI{0.00202}{\second}$, shortly after the peak. The figure illustrates how the wave propagates from the impact area on the rim across the surface of the machine casing. The results obtained with the FMM-ACA approach are similar and are therefore omitted here. A direct comparison is, however, given in Figure \ref{fig:disp_hammer}, where $u_{mag}$ is plotted at two distinct nodes for both approaches. 
Figure \ref{fig:disp_hammer1} shows $u_{mag}$ over time at node A, which belongs to an element on which input traction is prescribed. Figure \ref{fig:disp_hammer2} shows $u_{mag}$ at node B, which is located slightly further away from the impact, close to the underside of the casing. It can be clearly seen that the wave reaches node B slightly later, owing to its greater distance from the point of impact. Overall, both approaches are seen to agree very well at all times steps for both nodes A and B. More interesting, however, are the corresponding compression rates and the computing times of the two approaches. 

The $\mc{H}$-matrix ACA approach achieved a compression rate of $0.0137$ and required $80,291\,\si{\second}$, whereas the FMM-ACA achieved a compression rate of $0.0199$ and required $136,483\,\si{\second}$. Although the larger computation time of the FMM-ACA approach is consistent with the earlier results from subsection \ref{sec:cube}, the $\mc{H}$-matrix ACA approach now appears to achieve better compression rates, contrary to before. Overall, in order to obtain essentially the same results, $\mc{H}$-matrix ACA approach appears to be a better choice in terms of both compression and computation time. 

\subsubsection{Quadrupole excitation test}
\begin{figure}
    \centering
    \includegraphics[scale=1]{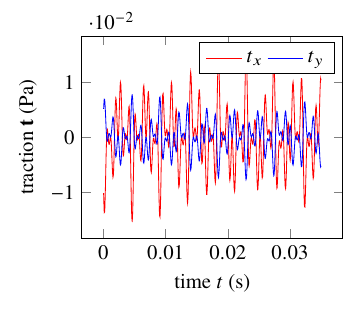}
    \caption{Harmonic radial traction}
    \label{fig:freq_traction}
\end{figure}

In the second test case, the machine is excited with a radial quadrupole mode. Again, a mixed boundary value problem is considered, with the base of the machine prescribed with homogeneous Dirichlet conditions. The remainder of the exterior casing is prescribed with zero traction. The Neumann boundary condition is imposed on the inner boundary of the casing in the form of a harmonic radial traction distribution with angular dependence $-2\varphi$, given by
\begin{equation}
\begin{aligned}
    t_{x}\kl{\x,t} &= 4\cdot 10^{-3} \rho \sum_{i=1}^{4}\sin{\kl{2\pi\omega_i t - 2\varphi}} \frac{x}{r} \\
    t_{y}\kl{\x,t} &= 7\cdot 10^{-3} \rho \sum_{i=1}^{4}\sin{\kl{2\pi\omega_i t - 2\varphi}} \frac{y}{r} \\
    t_{z}\kl{\x,t} &= 0 \\
    \textup{with } \varphi &= \tan^{-1}\kl{\frac{y}{x}}, \quad \omega_i \in \{247.7, 682.0, 900.42, 1335.86 \} [\textup{s}^{-1}],
\end{aligned}
\end{equation}
where the radial components are weighted by $x/r$ and $y/r$, respectively. Figure \ref{fig:freq_traction} shows the prescribed traction components $t_{x}\kl{\x,t}$ and $t_{y}\kl{\x,t}$ at one boundary element.

\begin{figure}[!h]
    \centering
    \begin{subfigure}[b]{0.49\textwidth}
        \includegraphics[scale = 0.2,
        trim={7cm 2cm 5cm 2cm},
        clip]{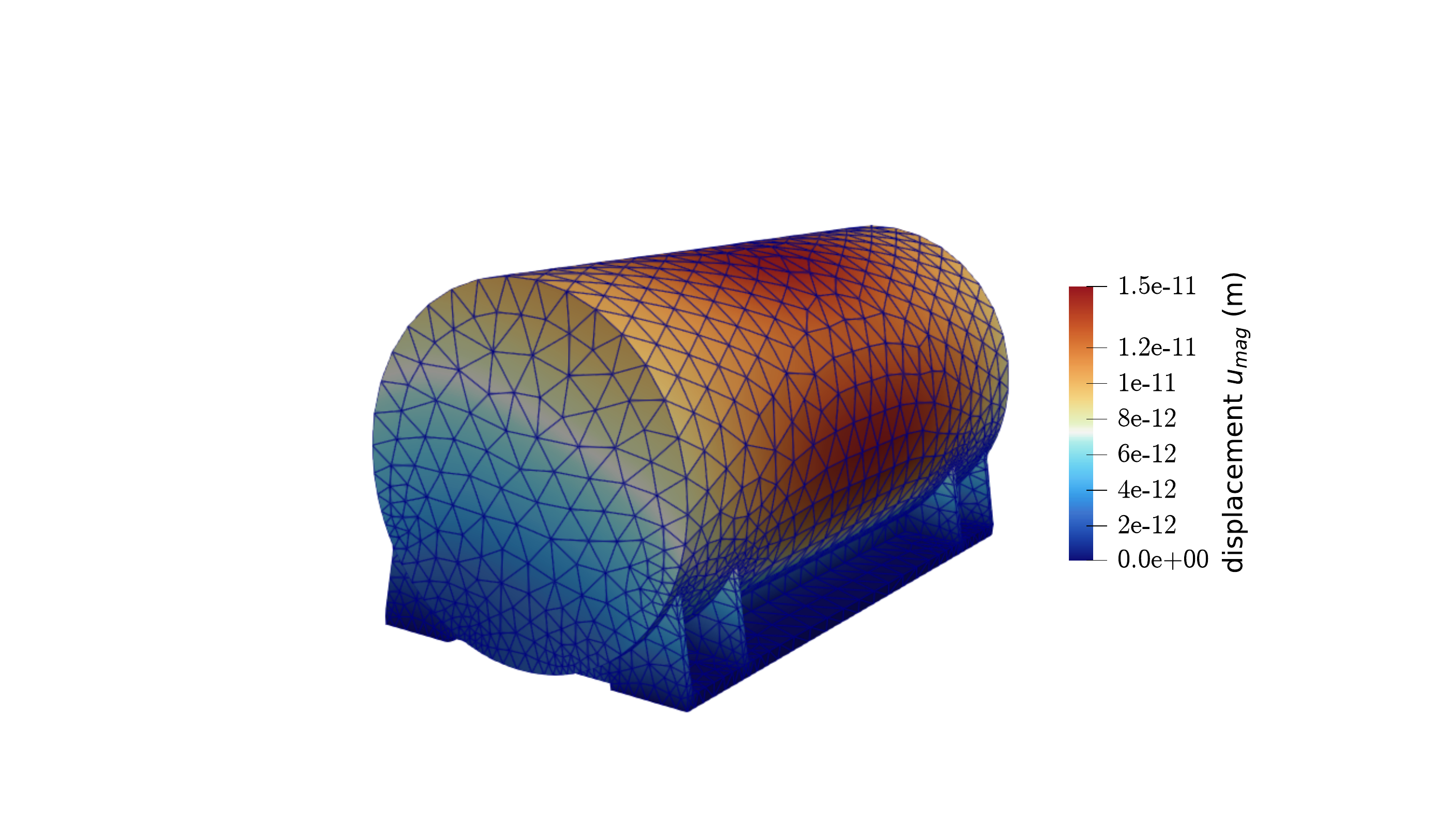}
        \caption{Time $t =  \SI{0.0100}{\second}$}
        \label{fig:freq1}
    \end{subfigure}
    \begin{subfigure}[b]{0.45\textwidth}
        \includegraphics[scale = 0.2,
        trim={7cm 2cm 5cm 2cm},
        clip]{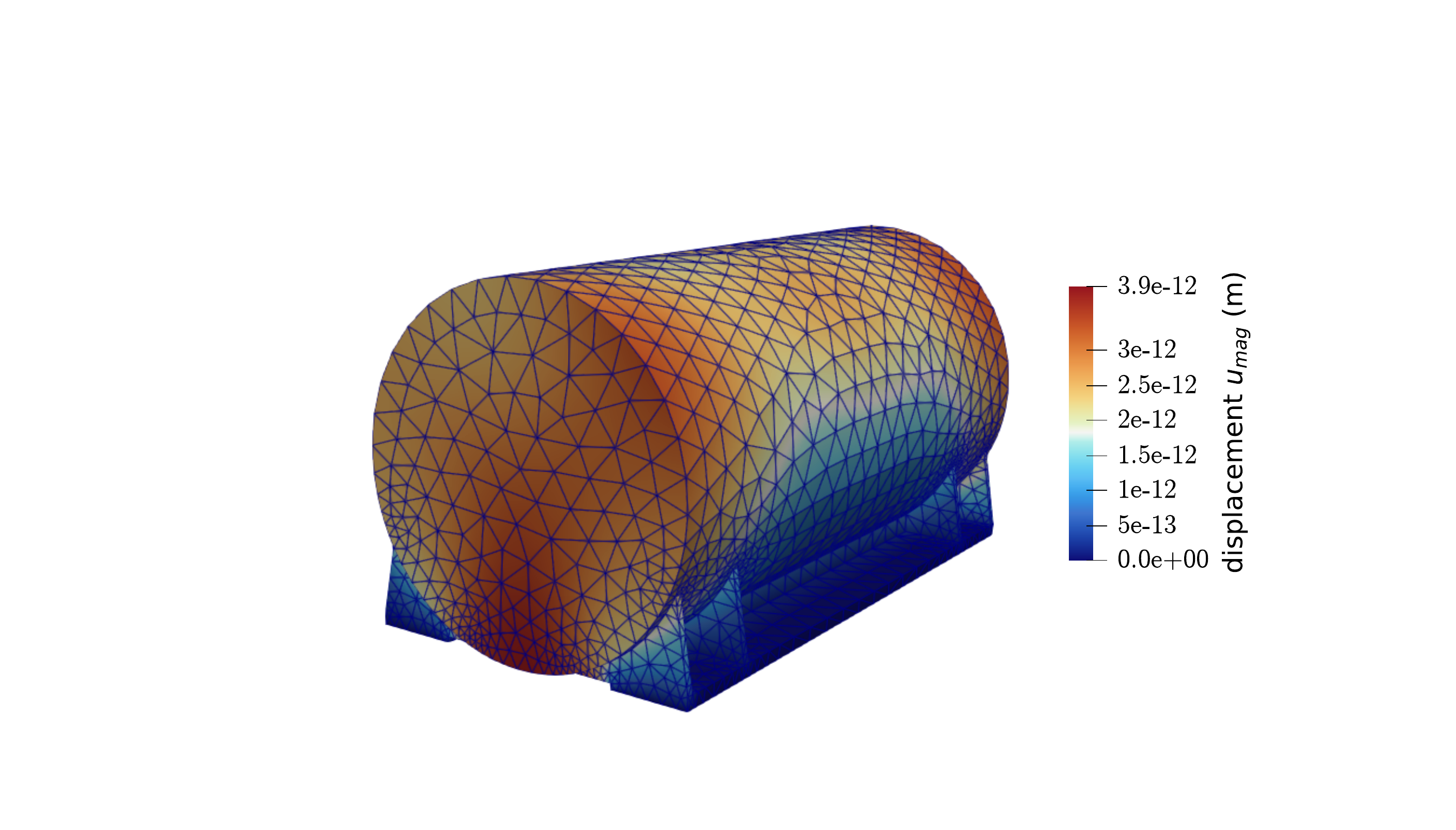}
        \caption{Time $t =  \SI{0.0276}{\second}$ }
        \label{fig:freq2}
    \end{subfigure}
    \caption{Magnitude of displacement $u_{mag}$ (in m) on the surface of the machine casing at different times (FMM-ACA)}
    \label{fig:machine_freq}
\end{figure}

The results obtained with the FMM-ACA approach are shown in Figure \ref{fig:machine_freq}.
Again, the displacement magnitude $u_{mag}$ is shown on the surface of the machine casing at two distinct time instants, namely at $t=  \SI{0.0100}{\second}$ and $t =  \SI{0.0276}{\second}$. At the first selected time instant, the vibration is concentrated mainly in the upper part of the casing, whereas the second time instant shows a more uniformly distributed wave pattern, as expected for a harmonic quadrupole excitation. The results obtained with the $\mc{H}$-matrix ACA approach are again very similar and are therefore omitted here. 

\begin{figure}[!h]
    \centering
    \begin{subfigure}[b]{0.45\textwidth}
        \includegraphics[scale = 1]{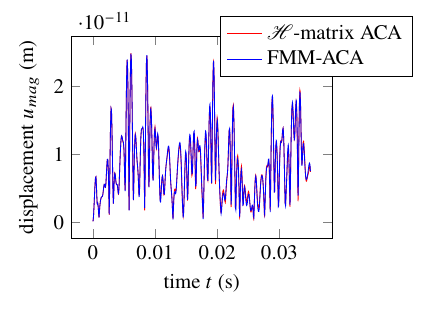}
        \caption{$u_{mag}$ at node A}
        \label{fig:disp_freq1}
    \end{subfigure}
    \begin{subfigure}[b]{0.45\textwidth}
        \includegraphics[scale = 1]{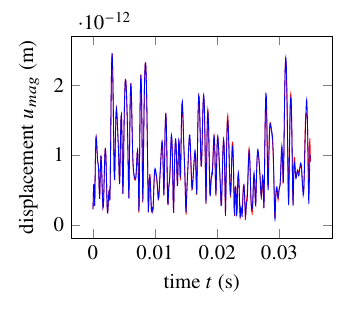}
        \caption{$u_{mag}$ at node B}
        \label{fig:disp_freq2}
    \end{subfigure}
    \caption{Magnitude of displacement $u_{mag}$ over time at two distinct nodes}
    \label{fig:disp_freq}
\end{figure}

Next, the displacement curves over time are presented in Figure \ref{fig:disp_freq}.
Similar to the previous case, the displacement magnitude $u_{mag}$ is plotted over time for both approaches, at two distinct nodes, A and B. In both Figures \ref{fig:disp_freq1} and \ref{fig:disp_freq2}, it can be clearly seen that FMM-ACA and $\mc{H}$-matrix ACA produce very similar results, in agreement with the previous examples. Although the two approaches select different frequencies from the gCQ in constructing the low-rank approximation, the chosen tolerance parameters are sufficient to yield nearly identical results. More relevant, however, is the comparison in terms of compression rates and total computing times. 

The $\mc{H}$-matrix ACA achieved a compression rate of $0.0033$ and required a total computing time of $486,508\,\si{\second}$. By contrast, the FMM-ACA achieved a compression rate of $0.0052$ and required $1,284,519.18\,\si{\second}$. As in the previous test case, the $\mc{H}$-matrix ACA approach appears to be a more favourable choice, since it provides better compression rates while requiring significantly less computing time to obtain similar results.

\section{Conclusions}
\label{sec6}
The multivariate adaptive cross approximation (3D-ACA) combined with either the $\mc{H}$-matrix based ACA or the Chebyshev-interpolation based FMM was presented for the acceleration of a gCQ-based time-domain BEM for three-dimensional elastodynamics. While ACA and FMM provide data-sparse approximations in the spatial dimension, the 3D-ACA yields a low-rank approximation of the resulting three-dimensional data array by adaptively reducing the number of complex frequencies required by the gCQ. In this way, compression is achieved not only in space but also in the additional frequency dimension arising from the temporal discretization. The numerical examples showed that both proposed approaches preserve the accuracy of the underlying formulation very well. For the academic unit-cube example, both the $\mc{H}$-matrix ACA and FMM-ACA approaches produced very similar convergence behaviour, confirming that the low-rank approximations in the spatial and frequency dimensions do not noticeably deteriorate the results. At the same time, both approaches achieved high compression rates, with the dominant contribution stemming from the 3D-ACA approximation of the three-dimensional data array.

A more detailed comparison of the two approaches revealed clear differences in their practical performance. For the academic example, the FMM-ACA approach resulted in slightly better compression rates, but required significantly longer computing times than the $\mc{H}$-matrix ACA approach. For the large-scale electric machine examples, the $\mc{H}$-matrix ACA approach proved to be more favourable overall, achieving both lower computing times and better compression rates while maintaining essentially the same solution quality. Overall, the proposed framework demonstrates that the additional frequency dimension introduced by gCQ can be treated efficiently by means of low-rank approximation, thereby making large-scale time-domain BEM computations in three-dimensional elastodynamics feasible.

\bigskip
\noindent\textbf{Acknowledgement} \hspace*{0.5em}This work is supported by the joint DFG/FWF Collaborative Research Centre CREATOR (DFG: Project-ID 492661287/TRR 361; FWF: 10.55776/F90) at TU Darmstadt, TU Graz and JKU Linz.


\appendix

\section{Analytical elastodynamic solution}\label{apndxA}
\setcounter{equation}{0}
\numberwithin{equation}{section}

The analytic full space solutions for displacement and traction fields are constructed according to \cite{eringen74}. The solution describes the vector-valued wave propagation of a unit impulse originating from a source point $\mathbf{p}$. The displacement field $\ubf\kl{\x,t}$ is first given by 
\begin{equation}
    u_i\kl{\x,t} = U'_{ij}\kl{\mathbf{r},t}d_j,
\end{equation}
where $\mathbf{r}=\mathbf{p}-\x$ and $\mathbf{d}$ is the direction vector. The function $U'_{ij}\kl{\mathbf{r},t}$ is a second-order displacement tensor given by
\begin{equation}\label{eq:apndxU}
  U'_{ij}\kl{\mathbf{r},t} = \frac{1}{4 \pi \rho} \kl{\overset{0}{f_{ij}}\kl{\mathbf{r}} \int_{1/c_1}^{1/c_2} \lambda f\kl{t-\lambda r } \mathrm{d}\lambda + \overset{1}{f_{ij}}\kl{\mathbf{r}} f\kl{t - \frac{r}{c_1}} +  \overset{2}{f_{ij}}\kl{\mathbf{r}} f\kl{t - \frac{r}{c_2}} }, 
\end{equation}
with 
\begin{equation*}
   \overset{0}{f_{ij}}\kl{\mathbf{r}} = \frac{1}{r} \kl{3r_ir_j-\delta_{ij}}, \quad \overset{1}{f_{ij}}\kl{\mathbf{r}} = \frac{r_i r_j}{rc_1^2}, \quad \overset{2}{f_{ij}}\kl{\mathbf{r}} = \frac{1}{r c_2^2}\kl{\delta_{ij} - r_ir_j}.  
\end{equation*}
The integral in \eqref{eq:apndxU} is evaluated as 
\begin{equation}
\begin{aligned}
   \int_{1/c_1}^{1/c_2} \lambda f\kl{t-\lambda r } \mathrm{d}\lambda &= \frac{1}{2ar^2} \bigg( \textup{exp} \kl{-q_1^2\kl{r}} - \textup{exp}\kl{-q_2^2\kl{r}}    \\
     &+ \sqrt{a\pi} \kl{ab-t} \kl{\textup{erf}\kl{q_1\kl{r}} - \textup{erf}\kl{q_2\kl{r}}   }  \bigg), 
\end{aligned}
\end{equation}
where 
\begin{equation}
    \textup{erf}\kl{x} = \frac{2}{\pi} \int_0^x \textup{exp}\kl{-t^2}\mathrm{d}t, \quad q_a\kl{r} = \frac{\sqrt{a}}{c_a}\kl{abc_a + r-c_at}.  
\end{equation}
The function $f\kl{r,t,c_a}$ is required to be at least twice differentiable with respect to the time variable $t$ and is chosen as
\begin{equation}
    f\kl{r,t,c_a} = \textup{exp}\kl{-a \kl{t - \frac{r}{c_a} - ab}^2 }.
\end{equation}
The full space traction field $\tbf\kl{\x,t}$ is then obtained as
\begin{equation}
    t_i\kl{\x,t} = -T'_{ijk}\kl{\mathbf{r},t} n_j\kl{\x}d_k,
\end{equation}
where $T'_{ijk}\kl{\mathbf{r},t}$ is a third-order tensor given by
\begin{equation}
    \begin{aligned}
        T'_{ijk}\kl{\mathbf{r},t} = \frac{\rho}{4 \pi} &\Bigg( \overset{0}{g}_{ijk}\kl{\mathbf{r}} \int_{1/c_1}^{1/c_2} \lambda f\kl{t-\lambda r} \mathrm{d}\lambda \\
        & + \overset{1}{g}_{ijk}\kl{\mathbf{r}} \kl{f\kl{t-\frac{r}{c_2}} - \kl{\frac{c_2}{c_1}}^2 f\kl{t-\frac{r}{c_1}}   } \\ 
        & + \overset{2}{g}_{ijk}\kl{\mathbf{r}} \kl{\dot{f}\kl{t-\frac{r}{c_2}} - \kl{\frac{c_2}{c_1}}^3 \dot{f}\kl{t-\frac{r}{c_1}}   } \\ 
        & + \overset{3}{g}_{ijk}\kl{\mathbf{r}} \kl{f\kl{t-\frac{r}{c_1}} - \frac{r}{c_1} \dot{f}\kl{t-\frac{r}{c_1}}   } \\ 
        & + \overset{4}{g}_{ijk}\kl{\mathbf{r}} \kl{f\kl{t-\frac{r}{c_2}} - \frac{r}{c_2} \dot{f}\kl{t-\frac{r}{c_2}}   }  \Bigg) \\
    \end{aligned}
\end{equation}
with
\begin{equation*}
    \begin{aligned}
        \overset{0}{g}_{ijk}\kl{\mathbf{r}} &= -\frac{6c_2^2}{r^2} \kl{5\frac{r_ir_jr_k}{r^3} - \frac{\delta_{ij}r_k + \delta_{ik}r_j + \delta_{jk} r_i}{r}   } \\ 
        \overset{1}{g}_{ijk}\kl{\mathbf{r}} &= \frac{2}{r^2} \kl{6\frac{r_ir_jr_k}{r^3} - \frac{\delta_{ij}r_k + \delta_{ik}r_j + \delta_{jk} r_i}{r}   } \\ 
        \overset{2}{g}_{ijk}\kl{\mathbf{r}} &= \frac{2r_ir_jr_k}{r^4c_2} \\
        \overset{3}{g}_{ijk}\kl{\mathbf{r}} &= - \kl{1-2\kl{\frac{c_2}{c_1}^2}} \frac{\delta_{ij}r_k}{r^3} \\
        \overset{4}{g}_{ijk}\kl{\mathbf{r}} &= -\frac{1}{r^2} \kl{\frac{\delta_{ik}r_j}{r} + \frac{\delta_{jk}r_i}{r}  }.
    \end{aligned}
\end{equation*}

For the unit-cube example, the parameters are chosen as $a = 10$ and $b = -0.055$. 

%% file: PreprintArxiv.bbl
\begin{thebibliography}{10}

\bibitem{sayas16}
F.-J. Sayas, {\em Retarded Potentials and Time Domain Boundary Integral
  Equations}, vol.~50 of {\em Springer Series in Computational Mathematics}.
\newblock Cham: Springer, 2016.
\newblock A road map.

\bibitem{costabel04}
M.~Costabel, {\em Time-Dependent Problems with the Boundary Integral Equation
  Method}, ch.~25, pp.~703--721.
\newblock John Wiley \& Sons, Ltd, 2004.

\bibitem{cruse68}
T.~Cruse and F.~Rizzo, ``A direct formulation and numerical solution of the
  general transient elastodynamic problem. {I},'' {\em Journal of Mathematical
  Analysis and Applications}, vol.~22, no.~1, pp.~244--259, 1968.

\bibitem{dominguez78}
J.~Dom{\'i}nguez and J.~M. Roesset, ``Dynamic stiffness of rectangular
  foundations,'' Tech. Rep. R78-20, Massachusetts Institute of Technology,
  Department of Civil Engineering, Cambridge, MA, Aug. 1978.

\bibitem{mansur83}
W.~J. Mansur, {\em A Time-Stepping Technique to Solve Wave Propagation Problems
  Using the Boundary Element Method}.
\newblock {P}h.{D}. {T}hesis, University of Southampton, 1983.

\bibitem{antes85}
H.~Antes, ``A boundary element procedure for transient wave propagations in
  two-dimensional isotropic elastic media,'' {\em Finite Elements in Analysis
  and Design}, vol.~1, no.~4, pp.~313--322, 1985.

\bibitem{israil90}
A.~Israil and P.~Banerjee, ``Two-dimensional transient wave-propagation
  problems by time-domain {BEM},'' {\em International Journal of Solids and
  Structures}, vol.~26, no.~8, pp.~851--864, 1990.

\bibitem{lubich1}
C.~Lubich, ``Convolution quadrature and discretized operational calculus.
  {I},'' {\em Numerische Mathematik}, vol.~52, no.~2, pp.~129--145, 1988.

\bibitem{lubich2}
C.~Lubich, ``Convolution quadrature and discretized operational calculus.
  {II},'' {\em Numerische Mathematik}, vol.~52, no.~4, pp.~413--425, 1988.

\bibitem{schanz97}
M.~Schanz and H.~Antes, ``A new visco- and elastodynamic time domain boundary
  element formulation,'' {\em Computational Mechanics}, vol.~20, no.~5,
  pp.~452--459, 1997.

\bibitem{banjai12}
L.~Banjai and M.~Schanz, ``Wave propagation problems treated with convolution
  quadrature and {BEM},'' in {\em Fast Boundary Element Methods in Engineering
  and Industrial Applications}, vol.~63 of {\em Lecture Notes in Applied and
  Computational Mechanics}, pp.~145--184, Heidelberg: Springer, 2012.

\bibitem{lopez13}
M.~L{\'o}pez-Fern{\'a}ndez and S.~Sauter, ``Generalized convolution quadrature
  with variable time stepping,'' {\em IMA Journal of Numerical Analysis},
  vol.~33, no.~4, pp.~1156--1175, 2013.

\bibitem{lopez15}
M.~L{\'o}pez-Fern{\'a}ndez and S.~Sauter, ``Generalized convolution quadrature
  with variable time stepping. part {II}: Algorithm and numerical results,''
  {\em Applied Numerical Mathematics}, vol.~94, pp.~88--105, 2015.

\bibitem{lopez16}
M.~L{\'o}pez-Fern{\'a}ndez and S.~Sauter, ``Generalized convolution quadrature
  based on {Runge--Kutta} methods,'' {\em Numerische Mathematik}, vol.~133,
  no.~4, pp.~743--779, 2016.

\bibitem{sauter17}
S.~A. Sauter and M.~Schanz, ``Convolution quadrature for the wave equation with
  impedance boundary conditions,'' {\em Journal of Computational Physics},
  vol.~334, pp.~442--459, 2017.

\bibitem{leitner21}
M.~Leitner and M.~Schanz, ``Generalized convolution quadrature based boundary
  element method for uncoupled thermoelasticity,'' {\em Mechanical Systems and
  Signal Processing}, vol.~150, p.~107234, 2021.

\bibitem{greengard87}
L.~Greengard and V.~Rokhlin, ``A fast algorithm for particle simulations,''
  {\em Journal of Computational Physics}, vol.~73, no.~2, pp.~325--348, 1987.

\bibitem{liu09}
Y.~Liu, {\em Fast Multipole Boundary Element Method}.
\newblock Cambridge: Cambridge University Press, 2009.
\newblock Theory and applications in engineering.

\bibitem{of05}
G.~Of, O.~Steinbach, and W.~L. Wendland, ``Applications of a fast multipole
  {Galerkin} boundary element method in linear elastostatics,'' {\em Computing
  and Visualization in Science}, vol.~8, no.~3--4, pp.~201--209, 2005.

\bibitem{phillips02}
J.~R. Phillips and J.~K. White, ``A precorrected-{FFT} method for electrostatic
  analysis of complicated 3-{D} structures,'' {\em IEEE Transactions on
  computer-aided design of integrated circuits and systems}, vol.~16, no.~10,
  pp.~1059--1072, 1997.

\bibitem{masters04}
N.~Masters and W.~Ye, ``Fast {BEM} solution for coupled 3{D} electrostatic and
  linear elastic problems,'' {\em Engineering Analysis with Boundary Elements},
  vol.~28, no.~9, pp.~1175--1186, 2004.

\bibitem{kapur97}
S.~Kapur and D.~E. Long, ``{IES\textsuperscript{3}}: A fast integral equation
  solver for efficient 3-dimensional extraction,'' in {\em Proceedings of the
  IEEE/ACM International Conference on Computer-Aided Design}, pp.~448--455,
  1997.

\bibitem{gope05}
D.~Gope and V.~Jandhyala, ``Efficient solution of {EFIE} via low-rank
  compression of multilevel predetermined interactions,'' {\em IEEE
  transactions on antennas and propagation}, vol.~53, no.~10, pp.~3324--3333,
  2005.

\bibitem{bebendorf03}
M.~Bebendorf and S.~Rjasanow, ``Adaptive low-rank approximation of collocation
  matrices,'' {\em Computing}, vol.~70, no.~1, pp.~1--24, 2003.

\bibitem{bebendorf08}
M.~Bebendorf, {\em Hierarchical Matrices}, vol.~63 of {\em Lecture Notes in
  Computational Science and Engineering}.
\newblock Berlin: Springer-Verlag, 2008.
\newblock A means to efficiently solve elliptic boundary value problems.

\bibitem{haider19}
A.~M. Haider and M.~Schanz, ``Adaptive cross approximation for {BEM} in
  elasticity,'' {\em Journal of Theoretical and Computational Acoustics},
  vol.~27, no.~1, 2019.
\newblock Paper No. 1850060, 19 pages.

\bibitem{messner10}
M.~Messner and M.~Schanz, ``An accelerated symmetric time-domain boundary
  element formulation for elasticity,'' {\em Engineering Analysis with Boundary
  Elements}, vol.~34, no.~11, pp.~944--955, 2010.

\bibitem{maruyama16}
T.~Maruyama, T.~Saitoh, T.~Q. Bui, and S.~Hirose, ``Transient elastic wave
  analysis of 3-{D} large-scale cavities by fast multipole {BEM} using implicit
  {Runge--Kutta} convolution quadrature,'' {\em Computer Methods in Applied
  Mechanics and Engineering}, vol.~303, pp.~231--259, 2016.

\bibitem{takahashi04}
T.~Takahashi, N.~Nishimura, and S.~Kobayashi, ``A fast {BIEM} for
  three-dimensional elastodynamics in time domain,'' {\em Engineering Analysis
  with Boundary Elements}, vol.~28, no.~2, pp.~165--180, 2004.

\bibitem{ergin98}
A.~A. Ergin, B.~Shanker, and E.~Michielssen, ``Fast evaluation of
  three-dimensional transient wave fields using diagonal translation
  operators,'' {\em Journal of Computational Physics}, vol.~146, no.~1,
  pp.~157--180, 1998.

\bibitem{bebendorf11}
M.~Bebendorf, ``Adaptive cross approximation of multivariate functions,'' {\em
  Constructive Approximation}, vol.~34, no.~2, pp.~149--179, 2011.

\bibitem{bebendorf13}
M.~Bebendorf, A.~K{\"u}hnemund, and S.~Rjasanow, ``An equi-directional
  generalization of adaptive cross approximation for higher-order tensors,''
  {\em Applied Numerical Mathematics}, vol.~74, pp.~1--16, 2013.

\bibitem{fong09}
W.~Fong and E.~Darve, ``The black-box fast multipole method,'' {\em Journal of
  Computational Physics}, vol.~228, no.~23, pp.~8712--8725, 2009.

\bibitem{schanz26}
M.~Schanz, V.~L. Keshava, and H.~De~Gersem, ``Comparison of
  {$\mathcal{H}$}-matrix- and {FMM}-based {3D-ACA} for a time-domain boundary
  element method,'' {\em Computational Mechanics}, 2026.

\bibitem{rjasanow17}
S.~Rjasanow and L.~Weggler, ``Matrix valued adaptive cross approximation,''
  {\em Mathematical Methods in the Applied Sciences}, vol.~40, no.~7,
  pp.~2522--2531, 2017.

\bibitem{chaillat17}
S.~Chaillat, L.~Desiderio, and P.~Ciarlet, ``Theory and implementation of
  {$\mathcal{H}$}-matrix based iterative and direct solvers for {Helmholtz} and
  elastodynamic oscillatory kernels,'' {\em Journal of Computational Physics},
  vol.~351, pp.~165--186, 2017.

\bibitem{seibel22}
D.~Seibel, ``Boundary element methods for the wave equation based on
  hierarchical matrices and adaptive cross approximation,'' {\em Numerische
  Mathematik}, vol.~150, no.~2, pp.~629--670, 2022.

\bibitem{dominguez93}
J.~Dom{\'i}nguez, {\em Boundary Elements in Dynamics}.
\newblock International Series on Computational Engineering, Southampton:
  Computational Mechanics Publications, 1993.

\bibitem{bonnet99}
M.~Bonnet, {\em Boundary Integral Equation Methods for Solids and Fluids}.
\newblock Chichester: John Wiley \& Sons, 1999.

\bibitem{graffi54}
D.~Graffi, ``{\"U}ber den {R}eziprozit{\"a}tssatz in der {D}ynamik der
  elastischen {K}{\"o}rper,'' {\em Ingenieur-Archiv}, vol.~22, no.~1,
  pp.~45--46, 1954.

\bibitem{wheeler68}
L.~T. Wheeler and E.~Sternberg, ``Some theorems in classical elastodynamics,''
  {\em Archive for Rational Mechanics and Analysis}, vol.~31, no.~1,
  pp.~51--90, 1968.

\bibitem{mantic93}
V.~Manti{\v{c}}, ``A new formula for the {C}-matrix in the {Somigliana}
  identity,'' {\em Journal of Elasticity}, vol.~33, no.~3, pp.~191--201, 1993.

\bibitem{hadamard23}
J.~Hadamard, {\em Lectures on {Cauchy's} Problem in Linear Partial Differential
  Equations}.
\newblock New Haven: Yale University Press, 1923.

\bibitem{kielhorn08}
L.~Kielhorn and M.~Schanz, ``Convolution quadrature method-based symmetric
  {Galerkin} boundary element method for 3-{D} elastodynamics,'' {\em
  International Journal for Numerical Methods in Engineering}, vol.~76, no.~11,
  pp.~1724--1746, 2008.

\bibitem{erichsen98}
S.~Erichsen and S.~A. Sauter, ``Efficient automatic quadrature in 3-{D}
  {Galerkin} {BEM},'' {\em Computer Methods in Applied Mechanics and
  Engineering}, vol.~157, no.~3--4, pp.~215--224, 1998.

\bibitem{duffy82}
M.~G. Duffy, ``Quadrature over a pyramid or cube of integrands with a
  singularity at a vertex,'' {\em SIAM Journal on Numerical Analysis}, vol.~19,
  no.~6, pp.~1260--1262, 1982.

\bibitem{schanz01}
M.~Schanz, {\em Wave Propagation in Viscoelastic and Poroelastic Continua: A
  Boundary Element Approach}, vol.~2 of {\em Lecture Notes in Applied
  Mechanics}.
\newblock Berlin: Springer, 2001.

\bibitem{eringen74}
A.~C. Eringen and E.~S. Suhubi, {\em Elastodynamics, Volume {I}: Finite
  Motions}.
\newblock New York: Academic Press, 1974.

\bibitem{kramer26}
T.~Kramer, B.~Marussig, and M.~Schanz, ``A higher-order time-domain boundary
  element formulation based on isogeometric analysis and the convolution
  quadrature method,'' {\em Computer Methods in Applied Mechanics and
  Engineering}, vol.~451, 2026.
\newblock Paper No. 118609, 18 pages.

\bibitem{weilharter11}
B.~Weilharter, O.~Biro, H.~Lang, and S.~Rainer, ``Computation of the noise
  radiation of an induction machine using {3D} {FEM}/{BEM},'' {\em COMPEL --
  The International Journal for Computation and Mathematics in Electrical and
  Electronic Engineering}, vol.~30, no.~6, pp.~1737--1750, 2011.

\end{thebibliography}
